# CONDITIONAL LIMIT THEOREMS FOR REGULATED FRACTIONAL BROWNIAN MOTION

By Hernan Awad and Peter Glynn

*University of Miami and Stanford University*

We consider a stationary fluid queue with fractional Brownian motion input. Conditional on the workload at time zero being greater than a large value $b$, we provide the limiting distribution for the amount of time that the workload process spends above level $b$ over the busy cycle straddling the origin, as $b \to \infty$. Our results can be interpreted as showing that long delays occur in large *clumps* of size of order $b^{2-1/H}$. The conditional limit result involves a finer scaling of the queueing process than fluid analysis, thereby departing from previous related literature.

**1. Introduction.** In the past ten years, there has been great interest in analyzing the performance of queues when the incoming traffic exhibits long-range dependence and self-similarity. This trend was motivated by measurements and statistical analysis of traffic in communication networks: after the initial findings of Leland et al. [23] in Ethernet traffic, numerous studies have shown that long-range dependent traffic is ubiquitous in high-speed communication networks, and have offered partial explanations for the origin of this phenomenon (see e.g., Crovella and Bestavros [7]).

One model that has received significant attention is that of a fluid queue that receives fractional Brownian motion (fBM) input—the so-called *fractional Brownian storage* (Norros [30]). As a traffic model, fBM is attractive because it is a stylized model (i.e., a low parameter tractable model) that is widely believed to be representative of long-range dependent, light-tailed (LRD-LT) traffic; see e.g., Norros [31] and Erramilli, Narayan and Willinger [11]. The tractability comes from the Gaussian self-similar characteristics of fBM, while its ability to approximate LRD-LT traffic is supported by various limit theorems, the majority of which consider superpositions of *on-off* fluid sources with heavy-tailed connection length distributions; see Konstantopoulos and Lin [19], Willinger et al. [39], Taqqu, Willinger and Sherman [35], Heath, Resnick and Samorodnitsky [14], Mikosch et al. [27] and Whitt [37, 38]. (Possible explanations for the origin of heavy-tailed distributions in communication networks are given by Fiorini, Sheahan and Lipsky [12] and Jelenkovic and Tan [17].) In this paper, we focus on an infinite-buffer fluid queue fed by fBM input.

We are interested in the "local" behavior of the queue in the neighborhood of a "typical" time at which a customer (packet) experiences a long delay. We would like to answer questions like the following: What is the conditional behavior of the traffic at such times? How many customers (packets) will experience large delays before congestion dissipates? This last question, about the size of the typical *clump* of long-delayed packets, provides a complement to the commonly used performance measure given by the steady-state probability of a packet experiencing long delays: the latter does not distinguish between having an occasional group of 10 packets suffering a long delay after an interval of $10^7$ packets served, versus having a group of $10^4$ close-by packets with long delays, out of every $10^{10}$ served. These situations can lead to drastically different "Quality of Service" (QoS) for some delay-sensitive applications (like audio traffic). The issue of *clumpiness* is also critical to the QoS in the finite-buffer counterpart of this model: applications like video traffic are robust with respect to the loss of isolated packets, but very sensitive to the loss of a large clump of nearby packets. Note that under first-come–first-served discipline, a customer (packet) experiencing a large delay is equivalent to the workload exceeding a large threshold. Hence, we state our results below by conditioning on a large exceedance for the workload process.

There are several results in the literature which are relevant to our work. Approximations for the steady-state exceedance probability for a queue fed by fBM were developed in a number of incremental steps: first a bound by Norros [30], then logarithmic-asymptotics by Duffield and O'Connell [9], sharper bounds by Massoulie and Simonian [26] and finally Narayan [28] obtained the exact asymptotics. The latter can also be obtained as a special case of a result by Hüsler and Piterbarg [15], who present exact asymptotics for the tail of the all-time maximum for a class of Gaussian processes with negative drift (note the steady-state workload is equal in distribution to the all-time maximum of negative drift fBM). The asymptotic behavior of the maximum workload over $[0, t]$ as $t \to \infty$ was derived by Zeevi and Glynn [40]. The conditional path (in fluid scale) that leads to a long busy-period was studied by Norros [32] and Mandjes, Mannersalo and Norros [24]. O'Connell and Procissi [34] and Chang, Yao and Zajic [4, 5, 6] give the conditional path (in fluid scale) that leads to an exceedance of a large level $b$. The typical fluid path turns out to be nonlinear; this is in contrast with the case of traffic that is short-range dependent and light-tailed (SRD-LT), where paths to a high exceedance are linear in great generality. More important for our purposes here, the typical fluid path to an exceedance of level $b$ (with $b$ large) is tangent at $b$, and does not spend any time above $b$. It follows that fluid analysis does not yield sufficient information about the distribution of the *clump* of packets that experience long delays; even the scaling behavior of the clump size is lost: the fluid limit only shows that the number of long-delayed packets must be of smaller order than $b$ itself. Hence, to address this problem, we need to look at the workload process on a finer time and space scale. This is one aspect in which our work departs from most of the existing analyses for the fBM-driven queue.

In the next section, we study a stationary version of the workload process, conditional on the workload at time 0 being greater than $b$. We derive a conditional limit result for the deviations of the workload from level $b$, under appropriate scaling, as $b \to \infty$. The relevant scaling consists in speeding up time by a factor of $b^{2-1/H}$ and compressing space by a factor of $b^{2H-1}$, where $H$ is the *Hurst parameter* (self-similarity index) of the fBM input. The limiting process is a fractional Brownian motion with symmetric negative polynomial drift (namely $-\kappa|t|^{2H}$), and started at a random level which is exponentially distributed; see Theorem 1. The result is proved under the assumption that $1/2 < H < 0.78$; we believe the same result holds for all $1/2 < H < 1$, but our proof technique does not cover all this range. (Recall fBM has LRD increments for $1/2 < H < 1$ and SRD increments for $0 < H < 1/2$.) To the best of our knowledge, this is the first conditional limit result for the fBM-driven queue that uses a finer-than-fluid scale. The highly non-Markovian structure of fBM (i.e., the process becomes Markov only when its entire history is incorporated into the state) makes this computation particularly challenging.

The conditional behavior described above is very different from the one observed in the case of SRD-LT traffic. Near the origin, the drift of the limiting process is close to zero; that is, conditional on a packet experiencing a large delay, the traffic's most natural tendency is not to restore itself immediately to its equilibrium behavior (as is the case with SRD-LT traffic). Thus, many other packets are likely to also experience large delays. This suggests one will observe big clumps. We formalize this by deriving the conditional limit law for the total sojourn above level $b$ during the busy-cycle straddling the origin (which is a proxy for the size of the clump of packets experiencing long delays). Consistent with the above scaling, the clump size scales up as $b^{2-1/H}$ as $b \to \infty$; see Corollary 2. Note that if one compares two fBM traffic sources, the one with higher Hurst parameter will exhibit asymptotically larger clump sizes, independent of the other parameters. Thus, in the words of Neidhardt and Wang [29], for this performance measure "Hurst is always naughty."

We believe the qualitative structure of the conditional limit theory derived here is likely to be inherited by LRD-LT traffic in general. We also believe that the scaling behavior and qualitative structure is representative of queues fed by LRD-LT traffic.

As mentioned earlier, in the finite-buffer counterpart of the model considered here clumpiness (of the loss process in this case) is also an important issue. There is a long tradition of approximating finite-buffer loss models by infinite-buffer models and associated exceedance computations. In particular, the exceedance probabilities results mentioned above are used to make buffer-size recommendations. In a companion paper, we consider a Brownian queue (i.e., the workload is modeled as regulated Brownian motion) and a GI/GI/1 with heavy-tails, and show that the qualitative structure of loss clumps is the same as that of the exceedance clumps. It is reasonable to expect a similar relation to be true here, suggesting that the loss clumps scale polynomially in the buffer size $b$ in the finite-buffer model.

From a traffic modeling standpoint, this paper contributes one building block towards the goal of developing a quantifiable measure of *burstiness* (which is of interest to the Internet traffic modeling community). Our view is that clumpiness of the loss process and burstiness of the incoming traffic are closely related; this paper's results therefore make rigorous the intuition that LRD-LT traffic is more *bursty* than SRD-LT traffic. In future work, we attempt related computations for stylized SRD-HT (short-range dependent, heavy-tailed) and LRD-HT (long-range dependent, heavy-tailed) traffic. These traffic types are also of interest in the network traffic modeling context: there is both empirical evidence and theoretical developments suggesting that stable-motions or fractional stable-motions can approximate high-speed network traffic in situations in which the arguments in favor of Gaussian models like fBM fail to apply—see, for example, Konstantopoulos and Lin [20, 21], Mikosch et al. [27], Tsoukatos and Makowski [36], and Laskin et al. [22].

**2. Model and main results.** We consider a single-server fluid queue, that receives fractional Brownian motion as input.

Let $Z = (Z(t): -\infty < t < \infty)$ be standard fractional Brownian motion (fBM) with Hurst parameter (self-similarity index) $H$. That is, $Z$ is a mean-zero Gaussian process with stationary increments and continuous sample paths, started at $Z(0) = 0$ and with covariance structure given by

$$\mathrm{E}Z(t)Z(s) = (|t|^{2H} + |s|^{2H} - |t-s|^{2H})/2.$$

The traffic process or *arrival process* to the queue, $A = (A(t): t \in \mathbb{R})$, is modeled as fBM with drift,

$$A(t) \triangleq \lambda t + \theta Z(t),$$

$t \in \mathbb{R}$. Here $\lambda$ corresponds to the mean rate of the traffic, and for $t > s$, $A(t) - A(s)$ represents the cumulative incoming traffic over the interval $(s,t]$. Let $c > \lambda$ denote the service rate, so that the so called *netput process* (or *free process*) is $X = (X(t): t \in \mathbb{R})$ given by

$$X(t) \triangleq A(t) - ct = -\mu t + \theta Z(t),$$

where $\mu \triangleq c - \lambda > 0$. The workload (or buffer content) process $W = (W(t): t \in \mathbb{R})$ is then obtained by applying the regulator mapping to $X$, which in this case translates to

$$W(t) \triangleq X(t) - \inf_{s \leq t} X(s),$$

$t \in \mathbb{R}$. Thus constructed, $W(t)$ represents the total amount of work present in the system at time $t$ (the *virtual waiting time*), and corresponds to stationary regulated fBM. Note that $W$ satisfies the relation

$$W(t) = X(t) - X(s) + W(s) \vee \left(X(s) - \inf_{s \leq u \leq t} X(u)\right), \qquad t > s.$$

We are interested in studying the "local" behavior of the workload process around a "typical" exceedance of level $b$. For this purpose, we define the process $Y^b$ given by

$$Y^b(t) \triangleq b^{-2H+1}(W(tb^{2-1/H}) - b),$$

$t \in \mathbb{R}$. Note that this is indeed a description of the local behavior near the origin: conditional on $\{W(0) > b\}$, the length of the busy-cycle straddling the origin scales up linearly in $b$, while the speed-up factor in $Y^b$ is $b^{2-1/H}$. Hence, for any fixed $t > 0$, as $b$ increases the path $(Y^b(s): -t \leq s \leq t)$ represents an ever smaller fraction of the busy-cycle straddling the origin.

Our main result gives the limiting law of $Y^b$ conditional on $W(0) > b$, as $b \to \infty$.

THEOREM 1. *If $1/2 < H < (\sqrt{17} - 1)/4$ then*

(1) $$(Y^b | W(0) > b) \Longrightarrow Y \qquad \text{as } b \to \infty,$$

*in $C(-\infty, \infty)$ with the topology of uniform convergence on compact sets, where $Y(t) = Y(0) - \kappa_1 |t|^{2H} + \theta Z(t)$, $t \in \mathbb{R}$, $Y(0) \perp\!\!\!\perp Z$, $Y(0) \stackrel{\mathcal{D}}{=} \frac{1}{\kappa_2} \mathrm{Exp}(1)$, $\kappa_1 = \frac{1}{2(1-H)}(\frac{\mu(1-H)}{H})^{2H}$, and $\kappa_2 = \frac{\mu^{2H}(1-H)^{2H-1}}{\theta^2 H^{2H}}$.*

As mentioned above, we believe the result is likely to hold for all $1/2 < H < 1$, but a portion of the argument in the proof only works for $1/2 < H < (\sqrt{17} - 1)/4$; see Step 12 in the proof of Theorem 1.

It is interesting to compare the case $H > 1/2$ with the case $H = 1/2$, i.e., with the result for the Brownian queue. There the limiting law is that of a Brownian motion with symmetric negative linear drift. Here, it is fractional Brownian motion with symmetric negative polynomial drift. The results are, to some degree, in agreement, since the conditional limiting law for the Brownian queue corresponds exactly to replacing formally $H = 1/2$ in Theorem 1. Qualitatively, however, there are significant differences. First, in the Brownian case the limiting law is obtained without the need to scale time and space. Also, in the Brownian case the drift is equal to $-\mu|t|$. In particular, the traffic process on $[0,\infty)$ obeys its usual (unconditional) law. In contrast, the polynomial drift $-\kappa_1 |t|^{2H}$ obtained in the LRD case has a derivative equal to zero at the origin. This, together with the scaling in place, indicates that from time 0 onwards, and for a period of order $b^{2-1/H}$, the traffic intensity will be close to 1. Thus, the traffic process does not revert immediately to its unconditioned dynamics (as in the Brownian case), but rather the traffic load remains higher than average, putting the queue in heavy-traffic, for a significant period.

The fact that the limiting process $Y$ has symmetric drift is also interesting, and perhaps somewhat surprising, since the typical path to level $b$ in fluid scale is not symmetric: the buildup from 0 to level $b$ is slower than the "draining" back to 0; see Chang, Yao and Zajic [5].

Next, we examine the effect on the size of the clump of long-delayed packets. The busy-cycle straddling the origin is $(\tau_L, \tau_R)$, where

$$\tau_R = \inf\{t \geq 0 : W(t) = 0\},$$
$$\tau_L = \sup\{t \leq 0 : W(t) = 0\}.$$

Let $V^b$ denote the total sojourn above level $b$ over $(\tau_L, \tau_R)$, i.e.,

$$V^b \triangleq \int_{\tau_L}^{\tau_R} 1_{(b,\infty)}(W(t))\, dt.$$

Its conditional limit distribution is given in the next result.

COROLLARY 2. *If $1/2 < H < (\sqrt{17} - 1)/4$, then*

$$(b^{-2+1/H} V^b | W(0) > b) \Longrightarrow V \triangleq \int_{-\infty}^{\infty} 1_{(0,\infty)}(Y(t))\, dt \qquad \text{as } b \to \infty,$$

*where $Y$ is as in* (1).

PROOF. Fix $a > 0$ and let $\varphi_a : C(-\infty, \infty) \to \mathbb{R}$ be given by $\varphi_a(y) = \int_{-a}^{a} 1_{(0,\infty)}(y(t))\, dt$. Put $G_a = \{y \in C(-\infty, \infty) : \lim_{t \to \infty} y(t) = \lim_{t \to -\infty} y(t) = -\infty, Leb(\{t \in (-a,a) : y(t) = 0\}) = 0\}$, where $Leb$ denotes Lebesgue measure. Note $\varphi_a(y_b) \to \varphi_a(y)$ whenever $y_b$ converges to $y \in G_a$ uniformly on compact sets and $P(Y \in G_a) = 1$; it then follows from Theorem 1 and the continuous mapping theorem that $(\varphi_a(Y^b)|W(0) > b) \Longrightarrow \varphi_a(Y)$ as $b \to \infty$.

Since $\varphi_a(Y) \nearrow V$ a.s. as $a \nearrow \infty$, it follows that $\lim_{a\to\infty} \lim_{b\to\infty} P(\varphi_a(Y^b) > x | W(0) > b) = P(V > x)$, $x \in \mathbb{R}$. On the other hand, $|P(\varphi_a(Y^b) > x | W(0) > b) - P(b^{-2+1/H}V^b > x | W(0) > b)| \leq P(F(a,b)|W(0) > b)$, where $F(a,b) = \{\exists a < t < \tau_R : Y^b(t) > 0\} \cup \{\exists \tau_L < t < -a : Y^b(t) > 0\} \cup \{|\tau_L| < ab^{2-1/H}\} \cup \{\tau_R < ab^{2-1/H}\}$. But $\lim_{a\to\infty} \lim_{b\to\infty} P(F(a,b)|W(0) > b) = 0$, whence the result follows. □

We see that the total sojourn above level $b$ scales as $b^{2-1/H}$. In the Brownian case, $H = 1/2$, the unscaled total sojourn above $b$ has a conditional limiting distribution. Hence, in a queue fed by fBM traffic with $H > 1/2$, larger delays tend to occur in larger clumps, whereas in the Brownian case the distribution of the typical exceedance clump of a given level is roughly independent of the magnitude of the said level.

It is interesting to note that the random variable $V$ plays a role in the prefactor of some asymptotic tail probabilities ($EV^{-1}$ being a variant of Pickands' constant), and also in other conditional limit results for sojourns above a high level of stochastic processes: For example, it appears in Theorems 3.3.1 and 5.5.1 of Berman [3], which deal, respectively, with high sojourns of stationary Gaussian processes and stationary diffusions over a finite interval; in both results, the probability that the (scaled) sojourn above $b$ is larger than $x$ converges to $E(1/\tilde{V}; \tilde{V} > x)$, where $\tilde{V}$ represents, like $V$ above, the total sojourn above zero of a Gaussian process with symmetric negative polynomial drift, started off at an exponentially distributed offset. (Those results are, however, quite different from that in Corollary 2: the workload process $W$ is neither Gaussian nor a diffusion, we cannot restrict it to a finite time interval, the limiting distribution of the scaled sojourn is directly the distribution of $V$ rather than a transformation of it, and in the case of $W$ the sojourn above $b$ scales up with $b$, whereas in Berman's results it scales down with $b$.)

An important tool in the proof of Theorem 1 is the following result, which states that if a negative-drift fBM is conditioned on hitting level $b$, the hitting time of $b$ will be "very near" its most likely value, $t^* = \frac{bH}{\mu(1-H)}$.

THEOREM 3. *Let $T_b \triangleq \inf\{t \geq 0 : X(t) \geq b\}$. For all $r > 0$ and $\varepsilon > 0$,*
$$P(|T_b - t^*|/b^{H+r} > \varepsilon | T_b < \infty) \longrightarrow 0$$
*as $b \to \infty$.*

For the proof, we refer the reader to Lemma 7 in Dieker [8], who proves a more general result; Theorem 3 corresponds to the particular case of fBM and (in his notation) $\delta(u) = O(u^{H+r})$.

Note that Theorem 3 is a refinement on the fluid behavior of $T_b$ obtained from large-deviations analysis, namely that $(|T_b - t^*|/b | T_b < \infty) \longrightarrow 0$ in probability as $b \to \infty$. It turns out that the fluid scale result is too coarse for our purposes. A further refinement can be found in Hüsler and Piterbarg [16], Theorem 1, who show that $((T_b - t^*)/b^H | T_b < \infty)$ converges weakly to a Gaussian random variable.

**3. Proofs.** Before proceeding with the proof of Theorem 1, we recall the notion of asymptotically equivalent events. We say that two collections of events $(E_1^b : b \geq 0)$ and $(E_2^b : b \geq 0)$ are *asymptotically equivalent* if

$$\mathrm{P}(E_1^b \triangle E_2^b) = o(\mathrm{P}(E_1^b))$$

as $b \to \infty$. In the proof, we will use the following well-known fact.

LEMMA 1 (Asmussen and Klüppelberg [2], Lemma 2.4(a)). *If $(E_1^b : b \geq 0)$ and $(E_2^b : b \geq 0)$ are asymptotically equivalent, then*

$$\|\mathrm{P}(\cdot | E_1^b) - \mathrm{P}(\cdot | E_2^b)\| \to 0$$

*as $b \to \infty$, where $\|\cdot\|$ denotes total variation distance.*

Throughout this section, we use the following notation: for any process $Z$ and $\sigma$-fields $\mathcal{F}$ and $\mathcal{G}$, we denote $Z^{\mathcal{F}}(t) = \mathrm{E}(Z(t)|\mathcal{F})$, $Z_{\mathcal{F}}(t) = Z(t) - Z^{\mathcal{F}}(t)$, and $Z_{\mathcal{G}}^{\mathcal{F}}(t) = Z^{\mathcal{F}}(t) - \mathrm{E}(Z^{\mathcal{F}}(t)|\mathcal{G})$. The complexity of the argument requires a significant amount of additional notation, which is introduced at the points when first needed; a notation summary table is included at the end for ease of reference.

PROOF OF THEOREM 1. To find the conditional limit law of $Y^b$ given $\{W(0) > b\}$, we study separately its intercept with the ordinate axis and its deviations from its value at the origin; that is, we study the random variable $Y^b(0) = b^{-2H+1}(W(0) - b)$—the scaled "overshoot" of $W$ over level $b$ at time 0—and the process $\hat{Y}^b = (\hat{Y}^b(t) : t \in \mathbb{R})$ given by

$$\hat{Y}^b(t) \triangleq b^{-2H+1} X(b^{2-1/H} t).$$

(Note $\hat{Y}^b(t) = Y^b(t) - Y^b(0)$ for $t \in (\tau_L, \tau_R)$.) We find the conditional limit laws of $\hat{Y}^b$ and $Y^b(0)$ separately, and then show that they are asymptotically independent as $b \to \infty$.

Observe that $\{W(0) > b\} = \{T^b < \infty\}$, where $T^b \triangleq \inf\{t \geq 0 : X(-t) < -b\}$. The event $\{T^b < \infty\}$ depends on the whole history of $X$ over $(-\infty, 0]$, and the law of $Y^b$ conditioned on such an event is hard to analyze. However, by Lemma 1, to prove that $(Y^b | T^b < \infty) \Longrightarrow Y$ it is enough to show that $(Y^b | A(b)) \Longrightarrow Y$, where $A(b)$ is some appropriately chosen asymptotically equivalent event for which the conditional law is more tractable.

To begin with, one can restrict attention to the history of $X$ over a finite interval. Put $t_L \triangleq t^* - b^{H+r}$ and $t_R \triangleq t^* + b^{H+r}$, for some small $r > 0$ (specified later). By Theorem 3, it follows that the event

$$A_0(b) \triangleq \{\exists s \in [-t_R, -t_L] : X(s) = -b\}$$

is asymptotically equivalent to $\{T^b < \infty\}$. The asymptotically equivalent event $A(b)$ that we construct is a subset of $A_0(b)$, and encompasses the intuition that, given $\{T^b < \infty\}$, $X$ will attain values close to $-b$ in a neighborhood of $-t^*$.

If one formally considers the law of $\hat{Y}^b$ conditioned on $X(-t^*) = -b$, then it is an easy exercise to verify that the conditional mean and covariance functions of this Gaussian process converge to those of the desired limit process $\hat{Y}$ as $b \to \infty$, where

$$\hat{Y}(t) = -\kappa_1 |t|^{2H} + \theta Z(t).$$

The first several steps in the proof show that this is still true when conditioning on $\mathcal{H} \triangleq \sigma((X(s): -t_R \leq s \leq -t_L))$, as long as $X$ remains close enough to $-b$ and satisfies some regularity conditions; making precise how close and the requisite regularity conditions will lead to the appropriate definition of $A(b)$.

The remainder of the argument is broken into "steps." A high-level overview can be obtained by reading the first portion of each, with the statement of the respective intermediate results, while skipping the details that follow. It is because of conditions imposed in Step 11 that the proof needs the assumption $H < (\sqrt{17} - 1)/4$; the intermediate results presented in the other steps are valid for $1/2 < H < 1$; see discussion in Step 12.

*Step* 1. *The conditional mean given the "endpoint"* $X(-t_L)$. Define the event

$$A_1(b) \triangleq \{X(-t_L) \in [-b, -b + b^\gamma]\},$$

$\widetilde{A}_1(b) = A_0(b) \cap A_1(b)$, and $\mathcal{K} \triangleq \sigma(X(-t_L))$. In this step we show that if

(2) $\quad\quad 0 < \gamma < H + r \quad \text{and} \quad 0 < r < (1-H)^2/H,$

then, given $M > 0$, there exists $\varepsilon_1 > 0$ such that, for all large enough $b$,

$$I_{\widetilde{A}_1(b)} |E(\hat{Y}^b(s)|\mathcal{K}) + \kappa_1 |t|^{2H}| \leq b^{-\varepsilon_1} I_{\widetilde{A}_1(b)},$$

$s \in (-M, M)$.

To prove this, put $\tilde{t}_L = b^{-2+1/H} t_L$ and note that

$E(\hat{Y}^b(s)|\mathcal{K})$

$$= b^{-2H+1} \left\{ (X(-t_L) - \mu t_L) \frac{\text{cov}(Z(-t_L), Z(sb^{2-1/H}))}{\text{var } Z(-t_L)} - \mu b^{2-1/H} s \right\}$$

$$= b^{-2H+1}\left\{\frac{1}{2}(X(-t_L) - \mu t_L)(|s/\tilde{t}_L|^{2H} + 1 - (1 + s/\tilde{t}_L)^{2H}) - \mu b^{2-1/H}s\right\}$$

$$= \frac{b^{-2H+1}}{2}(-b - \mu t^*)|s/\tilde{t}_L|^{2H}$$

$$+ b^{-2H+1}\left\{\frac{(-b - \mu t^*)}{2}(1 - (1 + s/\tilde{t}_L)^{2H}) - \mu\frac{s}{\tilde{t}_L}t^*\right.$$

(3)
$$\left. + \mu\frac{s}{\tilde{t}_L}(t^* - t_L) + \frac{1}{2}[X(-t_L) + b + \mu(t^* - t_L)]\right.$$

$$\left. \times (|s/\tilde{t}_L|^{2H} + 1 - (1 + s/\tilde{t}_L)^{2H})\right\}$$

$$= -\kappa_1|s|^{2H}(t^*/t_L)^{2H} - \frac{b^{2(1-H)}}{2(1-H)}(1 - (1 + s/\tilde{t}_L)^{2H} + 2Hs/\tilde{t}_L)$$

$$+ \mu\frac{s}{\tilde{t}_L}b^{1-H+r} + \frac{b^{1-2H}}{2}[X(-t_L) + b + \mu b^{H+r}]$$

$$\times (|s/\tilde{t}_L|^{2H} + 1 - (1 + s/\tilde{t}_L)^{2H}).$$

Since $|X(-t_L) + b| \leq b^\gamma < b^{H+r}$ on $\widetilde{A}_1(b)$, it follows that, on $\widetilde{A}_1(b)$,

$$|\mathrm{E}(\hat{Y}^b(s)|\mathcal{K}) + \kappa_1|s|^{2H}|$$

$$\leq \kappa_1|s|^{2H}[(1 + b^{H+r}/t_L)^{2H} - 1] + \frac{b^{2(1-H)}}{2(1-H)}|1 - (1 + s/\tilde{t}_L)^{2H} + 2Hs/\tilde{t}_L|$$

$$+ \frac{b^{1-H+r}}{2} \cdot \frac{|s|}{\tilde{t}_L}\left[2\mu + (1+\mu)\left|\frac{|s/\tilde{t}_L|^{2H} + 1 - (1 + s/\tilde{t}_L)^{2H}}{s/\tilde{t}_L}\right|\right]$$

$$\sim \kappa_1|s|^{2H} \cdot \frac{2Hb^{H+r}}{t_L} + \frac{b^{2(1-H)}}{2(1-H)}H(2H-1)(s/\tilde{t}_L)^2$$

$$+ \frac{b^{1-H+r}}{2} \cdot \frac{|s|}{\tilde{t}_L}[2\mu + 2H(1+\mu)]$$

$$= O(b^{r+H-1} + b^{-2(1-H)^2/H} + b^{-(1-H)^2/H+r}) = o(b^{-\varepsilon_1})$$

for any $0 < \varepsilon_1 < (1-H)^2/H - r$.

*Step 2. Conditional mean given both endpoints.* Define

$$\mathcal{G} \triangleq \sigma(X(-t_L), X(-t_R)),$$

$$A_2(b) = \{|Z_\mathcal{K}(-t_R)| \leq b^{\eta_3}\},$$

where $\eta_3 > 0$ is a constant (specified later), and $\widetilde{A}_2(b) = A_0(b) \cap A_1(b) \cap A_2(b)$. In this step, we show that if

(4) $\qquad \eta_3 < H^2 + (1-H)^2(1+H)/H + r(2H-1)$

then, given $M > 0$, there exists $\varepsilon_2 > 0$ such that, for all large enough $b$,

$$I_{\widetilde{A}_2(b)}|\mathrm{E}(\hat{Y}^b(s)|\mathcal{G}) - \mathrm{E}(\hat{Y}^b(s)|\mathcal{K})| \leq I_{\widetilde{A}_2(b)} b^{-\varepsilon_2},$$

$s \in (-M, M)$.

To verify this, put $D = 2b^{H+r}$ and observe that, on $A_2(b)$,

$$|\mathrm{E}(\hat{Y}^b(s)|\mathcal{G}) - \mathrm{E}(\hat{Y}^b(s)|\mathcal{K})|$$
$$= \frac{\theta b^{-2H+1}}{\operatorname{var} Z(-t_L) \operatorname{var} Z(-t_R) - \operatorname{cov}(Z(-t_L), Z(-t_R))^2} \cdot |Z_\mathcal{K}(-t_R)|$$
$$\times |\operatorname{var} Z(-t_L) \operatorname{cov}(Z(sb^{2-1/H}), Z(-t_R))$$
$$\quad - \operatorname{cov}(Z(-t_L), Z(sb^{2-1/H})) \operatorname{cov}(Z(-t_L), Z(-t_R))|$$
$$\leq \frac{\theta b^{-2H+1} b^{\eta_3}}{4(1+D/t_L)^{2H} - (1 + (1+D/t_L)^{2H} - (D/t_L)^{2H})^2}$$
$$\times \left| 2\left( \left|\frac{s}{\tilde{t}_L}\right|^{2H} + \left(1 + \frac{D}{t_L}\right)^{2H} - \left(1 + \frac{D}{t_L} + \frac{s}{\tilde{t}_L}\right)^{2H} \right) \right.$$
$$\left. - \left( \left|\frac{s}{\tilde{t}_L}\right|^{2H} + 1 - \left(1 + \frac{s}{\tilde{t}_L}\right)^{2H} \right) \left( 1 + \left(1 + \frac{D}{t_L}\right)^{2H} - \left(\frac{D}{t_L}\right)^{2H} \right) \right|$$
$$\sim \frac{\theta b^{\eta_3 - 2H+1}}{4(D/t_L)^{2H}} \cdot 4H(H-1) \cdot \frac{s}{\tilde{t}_L} \cdot \frac{D}{t_L}$$
$$= O(b^{\eta_3 - H^2 - (1-H)^2(1+H)/H - r(2H-1)})$$
$$= o(b^{-\varepsilon_2})$$

as $b \to \infty$, for any $0 < \varepsilon_2 < H^2 + (1-H)^2(1+H)/H + r(2H-1) - \eta_3$.

Step 3. *Conditional mean given $\mathcal{H}$.* Here, we show that the conditional mean of $\hat{Y}^b$ given $\mathcal{H}$ is close to its conditional mean given $\mathcal{G}$, as long as the values of $Z$ over the interval $[-t_R, -t_L]$ remain close enough to their conditional expectation given $\mathcal{G}$.

To be specific, define

$$A_3(b) = \{|Z_\mathcal{G}(-t_L - s)| \leq \max(s^H b^\delta, s^{H-\eta_1}), 0 \leq s \leq b^{H+r}\},$$
$$A_4(b) = \{|Z_\mathcal{G}(-t_R + s)| \leq \max(s^H b^\delta, s^{H-\eta_1}), 0 \leq s \leq b^{H+r}\},$$

where $\eta_1 > 0$ and $\delta > 0$ will be specified later, and put

$$A(b) = A_0(b) \cap A_1(b) \cap A_2(b) \cap A_3(b) \cap A_4(b).$$

In this step, we show that if

(5) $$\delta < (1-H)^2(1+H)/H - (1-H)(r+\eta_1) + r\eta_1,$$

then, given $M > 0$, there exists $\varepsilon_3 > 0$ such that, for all large enough $b$,
$$I_{A(b)}|\mathrm{E}(\hat{Y}^b(s)|\mathcal{H}) - \mathrm{E}(\hat{Y}^b(s)|\mathcal{G})| \leq I_{A(b)} b^{-\varepsilon_3}, \tag{6}$$
$s \in (-M, M)$.

(Note the term $s^{H-\eta_1}$ within the maximum controls the deviations of $Z$ from its conditional expectation given $\mathcal{G}$ in the immediate neighborhood of the endpoints $-t_R$ and $-t_L$; the term $s^H b^\delta$ controls the deviations of $Z$ from its conditional expectation given $\mathcal{G}$ away from the endpoints.)

To prove (6), we find a representation of the LHS in which the differences $Z(-t_L - s) - \mathrm{E}(Z(-t_L - s)|\mathcal{G})$ appear explicitly, in order to use the regularity conditions contained in $A_3(b)$ and $A_4(b)$. For this purpose, it is convenient to define the process $\widetilde{Z}$ by $\widetilde{Z}(u) \triangleq Z(-t_L + u) - Z(-t_L)$, $u \in \mathbb{R}$. Note $\widetilde{Z}$ is standard fBM, and that in terms of $\widetilde{Z}$, $\mathcal{H} = \sigma(\widetilde{Z}(t_L), (\widetilde{Z}(u): -D \leq u \leq 0))$ and $\mathcal{G} = \sigma(\widetilde{Z}(t_L), \widetilde{Z}(-D))$, where $D = 2b^{H+r}$ as before. Also,
$$\hat{Y}^b(u) = \theta b^{-2H+1}[\widetilde{Z}(t_L + b^{2-1/H}u) - \widetilde{Z}(t_L) - (\mu/\theta)b^{2-1/H}u].$$

Put $\mathcal{J} \triangleq \sigma((\widetilde{Z}(u): -D \leq u \leq 0))$. Then, for fixed $s, t \in (-M, M)$, and $b$ large enough that $t_L > b^{2-1/H}M$, we have the representation
$$\begin{pmatrix} \widetilde{Z}(t_L) \\ \widetilde{Z}(t_L + b^{2-1/H}s) \\ \widetilde{Z}(t_L + b^{2-1/H}t) \end{pmatrix} = \begin{pmatrix} \widetilde{Z}^{\mathcal{J}}(t_L) \\ \widetilde{Z}^{\mathcal{J}}(t_L + b^{2-1/H}s) \\ \widetilde{Z}^{\mathcal{J}}(t_L + b^{2-1/H}t) \end{pmatrix}
+ \left[\begin{pmatrix} a_{11} & 0 & 0 \\ a_{21} & a_{22} & 0 \\ a_{31} & a_{32} & a_{33} \end{pmatrix}\right] \begin{pmatrix} W_1 \\ W_2 \\ W_3 \end{pmatrix}, \tag{7}$$
where $(W_1, W_2, W_3)$ is i.i.d. standard Gaussian, $W_1 \perp\!\!\!\perp \mathcal{J}$, $(W_2, W_3) \perp\!\!\!\perp \mathcal{H}$ and the (deterministic) constants $(a_{ij})$ are as needed to match the covariances on both sides.

Hence,
$$\hat{Y}^b(s) = \theta b^{-2H+1}[\widetilde{Z}^{\mathcal{J}}(t_L + b^{2-1/H}s) - \widetilde{Z}^{\mathcal{J}}(t_L) \\
+ (a_{21}/a_{11} - 1)(\widetilde{Z}(t_L) - \widetilde{Z}^{\mathcal{J}}(t_L)) \\
+ a_{22}W_2 - (\mu/\theta)b^{2-1/H}s], \tag{8}$$
and it follows that
$$\mathrm{E}(\hat{Y}^b(s)|\mathcal{H}) - \mathrm{E}(\hat{Y}^b(s)|\mathcal{G}) \\
= \theta b^{-2H+1}[\widetilde{Z}^{\mathcal{J}}_{\mathcal{G}}(t_L + b^{2-1/H}s) - \widetilde{Z}^{\mathcal{J}}_{\mathcal{G}}(t_L) - (a_{21}/a_{11} - 1)\widetilde{Z}^{\mathcal{J}}_{\mathcal{G}}(t_L)]. \tag{9}$$

Gripenberg and Norros [13] show that, for $u \geq 0$, the conditional expectation $\widetilde{Z}^{\mathcal{J}}(u)$ can be expressed as
$$\widetilde{Z}^{\mathcal{J}}(u) = \int_{-D}^{0} g(D, u, v) \, d\widetilde{Z}(v),$$

where $g(D, u, -v) \triangleq C_H(D-v)^{-H+1/2}v^{-H+1/2}\int_0^u \frac{w^{H-1/2}(w+D)^{H-1/2}}{w+v}\,dw$ for $0 \le v \le D$ and $C_H = \sin(\pi(H-1/2))/\pi$. Here, the integral against fBM can be defined both as an $L^2$ limit or as an almost sure limit by approximating $g(D, u, \cdot)$ by a sequence of simple functions; for a discussion of integration against fBM see for example Duncan, Hu and Pasik-Duncan [10]. Moreover, a similar argument to that in Norros, Valkeila and Virtamo [33], Lemma 2.2, shows that an "integration by parts" formula holds, namely

$$
\begin{aligned}
\widetilde{Z}^{\mathcal{J}}(u) = &-\int_{-D}^0 g'(D, u, v)\left[\widetilde{Z}(v) + \frac{v}{D}\widetilde{Z}(-D)\right] dv \\
&- \frac{\widetilde{Z}(-D)}{D}\int_{-D}^0 g(D, u, v)\,dv \\
\triangleq &\lim_{\varepsilon \searrow 0}\left\{-\int_{-D+\varepsilon}^{-\varepsilon} g'(D, u, v)\left[\widetilde{Z}(v) + \frac{v}{D}\widetilde{Z}(-D)\right] dv \right. \\
&- \frac{\widetilde{Z}(-D)}{D}\int_{-D+\varepsilon}^{-\varepsilon} g(D, u, v)\,dv + g(D, u, -\varepsilon)\widetilde{Z}(-\varepsilon) \\
&- g(D, u, -D+\varepsilon)[\widetilde{Z}(-D+\varepsilon) - \widetilde{Z}(-D)] \\
&\left. - \varepsilon\frac{\widetilde{Z}(-D)}{D}[g(D, u, -\varepsilon) + g(D, u, -D+\varepsilon)]\right\},
\end{aligned}
\tag{10}
$$

where $g'(D, u, -v) \triangleq \partial g(D, u, -v)/\partial v$ and the limit is a.s. well defined by the Hölder continuity of fBM paths (cf. Norros, Valkeila and Virtamo [33], Lemma 2.2). Similarly, using a sequence of simple functions to approximate $g(D, u, \cdot)$ one can show that $\mathrm{E}(\widetilde{Z}^{\mathcal{J}}(u)|\mathcal{G})$ has a representation given by the RHS of (10) with $\widetilde{Z}(\cdot)$ replaced by $\widetilde{Z}^{\mathcal{G}}(\cdot)$. It then follows that

$$\widetilde{Z}^{\mathcal{J}}_{\mathcal{G}}(u) = \int_{-D}^0 g'(D, u, v)\widetilde{Z}_{\mathcal{G}}(v)\,dv \triangleq \lim_{\varepsilon \searrow 0}\int_{-D+\varepsilon}^{-\varepsilon} g'(D, u, v)\widetilde{Z}_{\mathcal{G}}(v)\,dv.$$

We use this representation of $\widetilde{Z}^{\mathcal{J}}_{\mathcal{G}}(u)$ to show that, if (5) holds, then on $A(b)$ the RHS in (9) is bounded by $O(b^{-\varepsilon_3})$, for some $\varepsilon_3 > 0$. The details are given in Lemma 2.

*Step* 4. *Convergence of the covariance function.* If $r < (1-H)^2/H$, then given $t, s \in \mathbb{R}$, there exists $\varepsilon_0 > 0$ such that, for all large enough $b$,

$$|\operatorname{cov}(\hat{Y}^b(t), \hat{Y}^b(s)|\mathcal{H}) - (\theta^2/2)(|t|^{2H} + |s|^{2H} - |t-s|^{2H})| \le b^{-\varepsilon_0}.$$

This is proved as Lemma 3 below, using the representation (7).

*Step* 5. *Convergence of finite-dimensional distributions of* $(\hat{Y}^b|A(b))$. Note that $A(b) \in \mathcal{H}$, and that $(\hat{Y}^b|\mathcal{H})$ is a Gaussian process. The convergence of the conditional mean and covariance functions from previous steps, together with an application of the Cramér–Wold device allow us to conclude that

$$(\hat{Y}^b|A(b)) \xrightarrow{\text{f.d.d.}} \hat{Y},$$

where

$$\hat{Y}(t) = -\kappa_1|t|^{2H} + \theta Z(t),$$

$t \in \mathbb{R}$, and $\xrightarrow{\text{f.d.d.}}$ denotes convergence of finite-dimensional distributions.

*Step* 6. *Tightness of* $(\hat{Y}^b|A(b))$. We now show that $\{(\hat{Y}^b|A(b)): b > 0\}$ is tight. We do so by showing that, given $M > 0$, there exists $C_M > 0$ such that, for all large enough $b$,

(11) $$I_{A(b)}\mathrm{E}((\hat{Y}^b(v) - \hat{Y}^b(u))^2|\mathcal{H}) \leq I_{A(b)}C_M|u-v|^{2H}$$

for all $u, v \in [-M, M]$, which implies

$$\mathrm{E}((\hat{Y}^b(v) - \hat{Y}^b(u))^2|A(b)) = \frac{\mathrm{E}I_{A(b)}\mathrm{E}((\hat{Y}^b(v) - \hat{Y}^b(u))^2|\mathcal{H})}{\mathrm{P}(A(b))} \leq C_M|u-v|^{2H},$$

which in turn is a sufficient condition for tightness with respect to the topology of uniform convergence on compact sets; see e.g., Karatzas and Shreve [18], page 64.

To verify (11), note that

$$\mathrm{E}((\hat{Y}^b(v) - \hat{Y}^b(u))^2|\mathcal{H}) = [E(\hat{Y}^b(v)|\mathcal{H}) - E(\hat{Y}^b(u)|\mathcal{H})]^2 \\ + \mathrm{var}(\hat{Y}^b(v) - \hat{Y}^b(u)|\mathcal{H}).$$

But, conditional on $\mathcal{H}$, $(\hat{Y}^b(v) - \hat{Y}^b(u))$ is Gaussian, and its conditional variance is bounded above by its (unconditional) variance, so that

$$\mathrm{var}(\hat{Y}^b(v) - \hat{Y}^b(u)|\mathcal{H}) \leq \theta^2|u-v|^{2H}.$$

Also, in Lemma 4, we use the representations (8) and (10) to show that the function $u \mapsto E(\hat{Y}^b(u)|\mathcal{H})$ is differentiable and that there exists a constant $\widetilde{C}_M$ such that, on $A(b)$,

$$\left|\frac{\partial E(\hat{Y}^b(u)|\mathcal{H})}{\partial u}\right| \leq \widetilde{C}_M,$$

$u \in [-M, M]$. In particular, on $A(b)$

$$|E(\hat{Y}^b(v)|\mathcal{H}) - E(\hat{Y}^b(u)|\mathcal{H})| \leq \widetilde{C}_M|u-v|$$

for $u, v \in [-M, M]$. Putting $C_M \triangleq \widetilde{C}_M^2 + 4M^2 + \theta^2$, (11), and hence, the tightness of $\{(\hat{Y}^b|A(b)): b > 0\}$, follows.

*Step* 7. *Limit law of* $\hat{Y}^b$. Given the tightness of $\{(\hat{Y}^b|A(b)): b > 0\}$, and the convergence of its finite-dimensional distributions, it follows that

$$(\hat{Y}^b|A(b)) \Longrightarrow \hat{Y},$$

where $\Longrightarrow$ denotes weak convergence in $C(-\infty,\infty)$ (with the topology of uniform convergence on compact sets).

We postpone to Step 11 the proof that $A(b)$ is asymptotically equivalent to $\{W(0) > b\}$. Assuming this for the moment, it follows by Lemma 1 that

(12) $$(\hat{Y}^b | W(0) > b) \Longrightarrow \hat{Y}.$$

*Step* 8. *Limit law of* $Y^b(0)$. We now show that

$$(Y^b(0) | W(0) > b) \Longrightarrow Y(0) \stackrel{\mathcal{D}}{=} \frac{1}{\kappa_2} \operatorname{Exp}(1).$$

The proof uses the exact asymptotic for the tail of $W(0)$ as given by Hüsler and Piterbarg [15], Corollary 2:

(13) $$\mathrm{P}(T^b < \infty) \sim K b^{(1-H)^2/H}(1 - \Phi(\kappa_3 b^{1-H})),$$

where $K$ is a constant (which they provide explicitly, in terms of Pickands' constant), $\kappa_3 = \sqrt{\kappa_2/(1-H)}$ and $\Phi(\cdot)$ is the c.d.f. of a standard Gaussian random variable. Since $(1 - \Phi(x)) \sim \phi(x)/x$ as $x \to \infty$, where $\phi(x) = \exp(-x^2/2)/\sqrt{2\pi}$ is the standard Gaussian density, we obtain

$$\mathrm{P}(T^b < \infty) \sim (K/\kappa_3) b^{-(1-H)(2-1/H)} \exp(-(\kappa_3^2/2) b^{2(1-H)}).$$

Hence,

$$\lim_{b \to \infty} \mathrm{P}(Y^b(0) > x | W(0) > b)$$
$$= \lim_{b \to \infty} \frac{\mathrm{P}(T^{b+xb^{2H-1}} < \infty)}{\mathrm{P}(T^b < \infty)}$$
$$= \lim_{b \to \infty} \left(\frac{b}{b + xb^{2H-1}}\right)^{(1-H)(2-1/H)}$$
$$\quad \times \exp\left(-\frac{\kappa_3^2 [(b + xb^{2H-1})^{2(1-H)} - b^{2(1-H)}]}{2}\right)$$
$$= \exp\left(-(\kappa_3^2/2) \lim_{b \to \infty} [(b + xb^{2H-1})^{2(1-H)} - b^{2(1-H)}]\right)$$
$$= \exp(-(\kappa_3^2/2) \cdot 2(1-H)x)$$
$$= \exp(-\kappa_2 x).$$

*Step* 9. *Asymptotic independence of* $Y^b(0)$ *and* $\hat{Y}^b$. We now show that $Y^b(0)$ and $\hat{Y}^b$ are asymptotically independent, in the sense that

$$\mathrm{P}(\hat{Y}^b \in \cdot, Y^b(0) > x | W(0) > b) \longrightarrow \mathrm{P}(\hat{Y} \in \cdot) \cdot \mathrm{P}((1/\kappa_2) \operatorname{Exp}(1) > x)$$

as $b \to \infty$.

To verify this, note

$$P(\hat{Y}^b \in \cdot, Y^b(0) > x | W(0) > b)$$
$$= P(\hat{Y}^b \in \cdot, W(0) > b + xb^{2H-1})/P(W(0) > b)$$
$$= P(\hat{Y}^b \in \cdot | W(0) > b + xb^{2H-1})P(W(0) > b + xb^{2H-1})/P(W(0) > b)$$
$$= P(\hat{Y}^b \in \cdot | W(0) > \tilde{b})P(Y^b(0) > x | W(0) > b),$$

where $\tilde{b} = b + xb^{2H-1}$. But we have shown above that

$$P(Y^b(0) > x | W(0) > b) \to P(Y(0) > x).$$

Hence, it is enough to show that $P(\hat{Y}^b \in \cdot | W(0) > \tilde{b}) \to P(\hat{Y} \in \cdot)$. Observe that $\hat{Y}^b = U^{\tilde{b}}$ where

$$U^{\tilde{b}}(t) = \left(\frac{b}{b + xb^{2H-1}}\right)^{2H-1} \hat{Y}^{\tilde{b}}(\beta(t,\tilde{b})),$$

$\beta(t,\tilde{b}) \triangleq t \cdot (\eta(\tilde{b})/\tilde{b})^{2-1/H}$ and $\eta$ is the inverse function of $b \mapsto b + xb^{2H-1}$. Note that (i) $\beta(\cdot,\tilde{b})$ converges to $t \mapsto t$ uniformly on compacts sets as $\tilde{b} \to \infty$, (ii) $(\frac{b}{b+xb^{2H-1}})^{2H-1} \nearrow 1$ as $\tilde{b} \to \infty$, (iii) $\tilde{b} \to \infty \iff b \to \infty$, and (iv) $\{(\hat{Y}^{\tilde{b}}|W(0) > \tilde{b}): \tilde{b} > b_0\}$ is tight (because of (12)). From (i)–(iv), it follows that for arbitrary $\varepsilon > 0$,

$$P(\rho(U^{\tilde{b}}, \hat{Y}^{\tilde{b}}) > \varepsilon | W(0) > b) \longrightarrow 0$$

as $b \to \infty$, where $\rho$ is the usual metric on $C_{(-\infty,\infty)}$. Since $(\hat{Y}^{\tilde{b}}|W(0) > \tilde{b}) \Longrightarrow \hat{Y}$ as $b \to \infty$, it then follows that $(U^{\tilde{b}}|W(0) > \tilde{b}) \Longrightarrow \hat{Y}$ as $b \to \infty$, i.e.,

$$(\hat{Y}^b | W(0) > \tilde{b}) \Longrightarrow \hat{Y},$$

as desired.

*Step* 10. *Limit law of* $Y^b$. It follows from the previous three steps that

$$((Y^b(0), \hat{Y}^b)|W(0) > b) \Longrightarrow (Y(0), \hat{Y})$$

in $\mathbb{R} \times C(-\infty,\infty)$, where $Y(0) \stackrel{\mathcal{D}}{=} (1/\kappa_2)\operatorname{Exp}(1)$, $Y(0) \perp\!\!\!\perp Z$. Since $Y^b(t) = Y^b(0) + \hat{Y}^b(t)$, $t \in (\tau_L, \tau_R)$, and since $-\tau_L, \tau_R \Longrightarrow \infty$ as $b \to \infty$, Theorem 1 then follows by the continuous mapping principle.

It only remains to be verified that $A(b)$ and $\{W(0) > b\}$ are asymptotically equivalent, as we claimed in Step 7 above. We verify this in the next step.

*Step* 11. *Asymptotic equivalence of* $A(b)$ *and* $\{T^b < \infty\}$. We need to show that $P\{A(b) \triangle \{T^b < \infty\}\} = o(P\{T^b < \infty\})$. By Theorem 3,

$$P(\{T^b < \infty\} \setminus A_0) = o(P(T^b < \infty)),$$

so it is enough to show that
$$P(A_0(b) \setminus A(b)) = o(P(T^b < \infty)).$$

In what follows (and for the rest of the paper) we frequently write $A_0, \ldots, A_4$, omitting the explicit dependence on $b$. Note that

$$A_0 \setminus A(b) \subset (A_0 \cap A_1^C) \cup (A_1 \cap A_2^C) \cup (A_1 \cap A_3^C) \cup (A_1 \cap A_4^C)$$
$$\subset A_5 \cup (A_0 \cap A_6) \cup (A_1 \cap A_2^C) \cup (A_1 \cap A_3^C) \cup (A_1 \cap A_4^C),$$

where $A_5 \triangleq \{X(-t_L) < -b\}$ and $A_6 \triangleq \{X(-t_L) > -b + b^\gamma\}$. Thus,

(14)
$$P(A_0 \setminus A(b)) \leq P(A_5) + P(A_0 \cap A_6) + P(A_1 \cap A_2^C)$$
$$+ P(A_1 \cap A_3^C) + P(A_1 \cap A_4^C).$$

We show that all the terms on the RHS in (14) are $o(P(T^b < \infty))$ as $b \to \infty$, as long as the parameters $\gamma$, $r$, $\eta_1$, $\eta_3$ and $\delta$ are chosen appropriately. We list the required conditions on the parameters below; the details of the arguments are relegated to Lemma 5.

That $P(A_5) = o(P(T^b < \infty))$ follows from Theorem 3.

The proof that $P(A_0 \cap A_6) = o(P(T^b < \infty))$ makes repeated use of Borell's inequality (Lemma 7) and requires that

(15) $$\gamma > H + r - (H - 1/2)(1 - H - r).$$

The proof that $P(A_1 \cap A_2^C) = o(P(T^b < \infty))$ uses the additional condition

(16) $$\eta_3 > H^2 + rH + (1 - H + r)/2.$$

The proof that $P(A_1 \cap A_3^C)$ and $P(A_1 \cap A_4^C)$ are $o(P(T^b < \infty))$ relies on Borell's inequality and the Hölder continuity of fractional Brownian paths, and uses the condition

(17) $$\delta > (1 - H + r)/2.$$

*Step* 12. *Final remarks.* If $1/2 < H < (\sqrt{17} - 1)/4$, then the parameters $r, \gamma, \eta_3, \delta$ and $\eta_1$ in the definition of $A(b)$ can be chosen so as to satisfy (2), (4), (5), (15)–(17).

It is because of conditions (16) and (17) that the proof does not work for all $1/2 < H < 1$: these conditions, when combined with (4) and (5), require that $(1-H)/2 < (1-H)^2(1+H)/H$, which in turn requires $H < (\sqrt{17}-1)/4$. The need for (16) and (17) may be a consequence of the style of proof rather than necessary conditions. If one is able to prove that $P(A_1 \cap A_3^C)$ and $P(A_2 \cap A_3^C)$ are $o(P(T^b < \infty))$ without imposing (16) and (17), then the proof of Theorem 1 would work for all $H \in (1/2, 1)$. $\square$

LEMMA 2. *Fix $s > 0$. If $\delta$ satisfies (5), then there exists $\varepsilon_3 > 0$ such that*

(18) $$I_{A(b)}|\mathrm{E}(\hat{Y}^b(s)|\mathcal{H}) - \mathrm{E}(\hat{Y}^b(s)|\mathcal{G})| \leq I_{A(b)} b^{-\varepsilon_3}$$

*for all large enough $b$.*

PROOF. We analyze separately each of the terms on the RHS of equation (9). Using that $g'(D, u, v) = g'(D/u, 1, v/u)/u$ (which is straightforward to verify), we note that on $A(b)$,

$$|\tilde{Z}_\mathcal{G}^\mathcal{J}(t_L)| = \left|\int_0^D g'(D, t_L, -v)\tilde{Z}_\mathcal{G}(-v)\,dv\right|$$

$$\leq \int_0^{D/t_L} |g'(D/t_L, 1, -x)\tilde{Z}_\mathcal{G}(-xt_L)|\,dx$$

$$\leq b^\delta t_L^H \left[\int_0^{D/2t_L} |g'(D/t_L, 1, -x)|x^{H-\eta_1}\,dx\right.$$

$$\left.+ \int_{D/2t_L}^{D/t_L} |g'(D/t_L, 1, -x)|(D/t_L - x)^{H-\eta_1}\,dx\right],$$

where the last step follows since, on $A(b)$, $|\tilde{Z}_\mathcal{G}(-v)| \leq \max(v^H b^\delta, v^{H-\eta_1})$ for $0 \leq v \leq D/2$ and $|\tilde{Z}_\mathcal{G}(-v)| \leq \max((D-v)^H b^\delta, (D-v)^{H-\eta_1})$ for $D/2 \leq v \leq D$. Since

$$g'(D, u, -v) = -(H - 1/2)g(D, u, -v)[v^{-1} - (D-v)^{-1}]$$

$$- C_H v^{-H+1/2}(D-v)^{-H+1/2}\int_0^u \frac{y^{H-1/2}(y+D)^{H-1/2}}{(y+v)^2}\,dy,$$

it follows that

$$|\tilde{Z}_\mathcal{G}^\mathcal{J}(t_L)|$$

$$\leq b^\delta t_L^H C_H$$

$$\times \left\{\int_0^{D/2t_L} x^{1/2-\eta_1}\left(\frac{D}{t_L} - x\right)^{-H+1/2}\right.$$

$$\times \int_0^1 \frac{y^{H-1/2}(y+D/t_L)^{H-1/2}}{(y+x)^2}\,dy\,dx$$

$$+ \int_{D/2t_L}^{D/t_L} x^{-H+1/2}\left(\frac{D}{t_L} - x\right)^{1/2-\eta_1}$$

$$\left.\times \int_0^1 \frac{y^{H-1/2}(y+D/t_L)^{H-1/2}}{(y+x)^2}\,dy\,dx\right.$$

$$+ \left(H - \frac{1}{2}\right) \int_0^{D/2t_L} x^{1/2-\eta_1} \left(\frac{D}{t_L} - x\right)^{-H+1/2}$$
$$\times \left[x^{-1} + \left(\frac{D}{t_L} - x\right)^{-1}\right] h_1\left(\frac{D}{t_L}, x\right) dx$$
$$+ \left(H - \frac{1}{2}\right) \int_{D/2t_L}^{D/t_L} x^{-H+1/2} \left(\frac{D}{t_L} - x\right)^{1/2-\eta_1}$$
$$\times \left[x^{-1} + \left(\frac{D}{t_L} - x\right)^{-1}\right] h_1\left(\frac{D}{t_L}, x\right) dx \bigg\},$$

where

$$h_1(D/t_L, x) \triangleq \int_0^1 \frac{y^{H-1/2}(y + D/t_L)^{H-1/2}}{y + x} dy \leq 4/(2H - 1).$$

Note also that

$$\int_0^1 \frac{y^{H-1/2}(y + D/t_L)^{H-1/2}}{(y + x)^2} dy \leq 4 x^{H-3/2} (D/t_L)^{H-1/2},$$

so that

$$|\widetilde{Z}_{\mathcal{G}}^{\mathcal{J}}(t_L)| \leq 4 b^\delta t_L{}^H C_H$$
$$\times \bigg\{ (D/t_L)^{H-1/2} \int_0^{D/2t_L} x^{H-1-\eta_1} (D/t_L - x)^{-H+1/2} dx$$
$$+ (D/t_L)^{H-1/2} \int_{D/2t_L}^{D/t_L} x^{-1} (D/t_L - x)^{1/2-\eta_1} dx$$
(19)
$$+ \int_0^{D/2t_L} x^{1/2-\eta_1} (D/t_L - x)^{-H+1/2} [x^{-1} + (D/t_L - x)^{-1}] dx$$
$$+ \int_{D/2t_L}^{D/t_L} x^{-H+1/2} (D/t_L - x)^{1/2-\eta_1} [x^{-1} + (D/t_L - x)^{-1}] dx \bigg\}.$$
$$\leq 4 b^\delta t_L{}^H C_H \cdot 15 (D/t_L)^{1-H-\eta_1}$$
$$\leq c_1 b^{\delta + H - (1-H-r)(1-H-\eta_1)}$$

for an appropriate constant $c_1$.

Similarly, on $A(b)$,

$$|\widetilde{Z}_{\mathcal{G}}^{\mathcal{J}}(t_L) - \widetilde{Z}_{\mathcal{G}}^{\mathcal{J}}(t_L - b^{2-1/H} s)|$$
$$= \left| \int_0^D [g'(D, t_L, -v) - g'(D, t_L - b^{2-1/H} s, -v)] \widetilde{Z}_{\mathcal{G}}(-v) \, dv \right|$$
$$\leq \int_0^{D/t_L} |[g'(D/t_L, 1, -x) - g'(D/t_L, 1 - s/\tilde{t}_L, -x)] \widetilde{Z}_{\mathcal{G}}(-x t_L)| \, dx$$

$$\leq b^\delta t_L^H \Bigg[ \int_0^{D/2t_L} |g'(D/t_L, 1, -x) - g'(D/t_L, 1 - s/\tilde{t}_L, -x)| x^{H-\eta_1} \, dx$$

$$+ \int_{D/2t_L}^{D/t^L} |g'(D/t_L, 1, -x)$$

$$- g'(D/t_L, 1 - s/\tilde{t}_L, -x)|(D/t_L - x)^{H-\eta_1} \, dx \Bigg],$$

where $\tilde{t}_L \triangleq b^{-2+1/H} t_L$. Rewriting,

$$|\widetilde{Z}_\mathcal{G}^\mathcal{J}(t_L) - \widetilde{Z}_\mathcal{G}^\mathcal{J}(t_L - b^{2-1/H} s)|$$

$$\leq b^\delta t_L^H C_H (H - 1/2)$$

$$\times \Bigg\{ \frac{2}{2H-1} \int_0^{D/2t_L} x^{1/2-\eta_1} \left(\frac{D}{t_L} - x\right)^{-H+1/2}$$

$$\times \int_{1 \wedge (1-s/\tilde{t}_L)}^{1 \vee (1-s/\tilde{t}_L)} \frac{y^{H-1/2}(y + D/t_L)^{H-1/2}}{(y+x)^2} \, dy \, dx$$

$$+ \frac{2}{2H-1} \int_{D/2t_L}^{D/t_L} x^{-H+1/2} \left(\frac{D}{t_L} - x\right)^{1/2-\eta_1}$$

$$\times \int_{1 \wedge (1-s/\tilde{t}_L)}^{1 \vee (1-s/\tilde{t}_L)} \frac{y^{H-1/2}(y + D/t_L)^{H-1/2}}{(y+x)^2} \, dy \, dx$$

$$+ \int_0^{D/2t_L} x^{1/2-\eta_1} \left(\frac{D}{t_L} - x\right)^{-H+1/2}$$

$$\times \left[ x^{-1} + \left(\frac{D}{t_L} - x\right)^{-1} \right] h_2(D/t_L, s/t_L, x) \, dx$$

$$+ \int_{D/2t_L}^{D/t_L} x^{-H+1/2} \left(\frac{D}{t_L} - x\right)^{1/2-\eta_1}$$

$$\times \left[ x^{-1} + \left(\frac{D}{t_L} - x\right)^{-1} \right] h_2(D/t_L, s/t_L, x) \, dx \Bigg\},$$

where

$$h_2(D/t_L, s/t_L, x) \triangleq \int_{1 \wedge (1-s/\tilde{t}_L)}^{1 \vee (1-s/\tilde{t}_L)} \frac{y^{H-1/2}(y + D/t_L)^{H-1/2}}{y+x} \, dy \leq 2|s|/\tilde{t}_L.$$

Since also

$$\int_{1 \wedge (1-s/\tilde{t}_L)}^{1 \vee (1-s/\tilde{t}_L)} \frac{y^{H-1/2}(y + D/t_L)^{H-1/2}}{(y+x)^2} \, dy \leq 2|s|/\tilde{t}_L,$$

it follows that

$$|\widetilde{Z}_\mathcal{G}^\mathcal{J}(t_L) - \widetilde{Z}_\mathcal{G}^\mathcal{J}(t_L - b^{2-1/H} s)|$$

$$\leq 2s\tilde{t}_L^{-1}b^\delta t_L{}^H C_H$$

$$\times \left\{ \int_0^{D/2t_L} x^{1/2-\eta_1}(D/t_L - x)^{-H+1/2}\, dx \right.$$

$$+ \int_{D/2t_L}^{D/t_L} x^{-H+1/2}(D/t_L - x)^{1/2-\eta_1}\, dx$$

(20)
$$+ \int_0^{D/2t_L} x^{1/2-\eta_1}(D/t_L - x)^{-H+1/2}[x^{-1} + (D/t_L - x)^{-1}]\, dx$$

$$\left. + \int_{D/2t_L}^{D/t_L} x^{-H+1/2}(D/t_L - x)^{1/2-\eta_1}[x^{-1} + (D/t_L - x)^{-1}]\, dx \right\}$$

$$\leq 2s\tilde{t}_L^{-1}b^\delta t_L{}^H C_H (D/t_L)^{1-H-\eta_1}$$

$$\leq c_2 b^{\delta + H - (1-H+r)(1-H-\eta_1) - (1/H-1)}$$

for an appropriate constant $c_2$ and all large enough $b$.

To analyze the factor $(a_{21}/a_{11} - 1)$ in the first term on the RHS in (9), note that

$$\frac{a_{21}}{a_{11}} = \frac{\mathrm{E}\widetilde{Z}(t_L)\widetilde{Z}(t_L + b^{2-1/H}s) - \mathrm{E}\widetilde{Z}^{\mathcal{J}}(t_L)\widetilde{Z}^{\mathcal{J}}(t_L + b^{2-1/H}s)}{\mathrm{E}\widetilde{Z}(t_L)^2 - \mathrm{E}\widetilde{Z}^{\mathcal{J}}(t_L)^2}.$$

The covariance $\mathrm{E}\widetilde{Z}^{\mathcal{J}}(a)\widetilde{Z}^{\mathcal{J}}(a+x)$ satisfies,

$$\mathrm{E}\widetilde{Z}^{\mathcal{J}}(a)\widetilde{Z}^{\mathcal{J}}(a+x) = a^{2H} f(D/a, 0, x)$$

for $|x| < a$, where

$$f(d, w, x) \triangleq H \int_0^d g(d, 1+w, -u)[(1+w+x+u)^{2H-1} - u^{2H-1}]\, du$$

(see Mannersalo [25]). It follows that

$$\left(\frac{a_{21}}{a_{11}} - 1\right) = \frac{v(s/\tilde{t}_L) - 1 + f(D/t_L, 0, 0) - f(D/t_L, 0, s/\tilde{t}_L)}{1 - f(D/t_L, 0, 0)},$$

where $v(x) \triangleq \frac{1}{2}(1 + (1+x)^{2H} - |x|^{2H})$, $|x| < 1$. Note that $(v(x) - 1)/x \to H$ as $x \to 0$ and it can be easily shown that $f(D/t_L, 0, 0) - f(D/t_L, 0, s/\tilde{t}_L) = o(s/\tilde{t}_L)$ as $b \to \infty$, whence there exists a constant $c_3$ such that

(21) $$\left|\frac{a_{21}}{a_{11}} - 1\right| \leq c_3 |s| b^{-(1/H-1)}.$$

We can then use the above bounds (19)–(21) on the RHS of (9) to conclude that there exists $c_4$ such that, on $A(b)$,

$$|\mathrm{E}(\hat{Y}^b(s)|\mathcal{H}) - \mathrm{E}(\hat{Y}^b(s)|\mathcal{G})| \leq c_4 b^{-2H+1} b^{\delta + H - (1-H+r)(1-H-\eta_1) - (1/H-1)}$$

$$= o(b^{-\varepsilon_3}),$$

where $0 < \varepsilon_3 < (1-H)^2(1+H)/H - \delta - (1-H)(r+\eta_1) + r\eta_1 > 0$.  □

LEMMA 3. *Fix $s, t > 0$. If $r < (1-H)^2/H$ then there exists $\varepsilon_2 > 0$ such that*

$$|\operatorname{cov}(\hat{Y}^b(t), \hat{Y}^b(s)|\mathcal{H}) - (\theta^2/2)(|t|^{2H} + |s|^{2H} - |t-s|^{2H})| \leq b^{-\varepsilon_2}$$

*for all large enough $b$.*

PROOF. Note from (7) that

$$\operatorname{cov}(\hat{Y}^b(s), \hat{Y}^b(t)|\mathcal{H})$$
$$= b^{2(1-2H)}\theta^2 a_{32}a_{22}$$
$$= b^{2(1-2H)}\theta^2 \mathrm{E}\widetilde{Z}(t_L + b^{2-1/H}s)\widetilde{Z}(t_L + b^{2-1/H}t)$$
$$\quad - b^{2(1-2H)}\theta^2 \mathrm{E}\widetilde{Z}^{\mathcal{J}}(t_L + b^{2-1/H}s)\widetilde{Z}^{\mathcal{J}}(t_L + b^{2-1/H}t)$$
$$\quad - b^{2(1-2H)}\theta^2 \cdot \frac{(\mathrm{E}\widetilde{Z}(t_L)\widetilde{Z}(t_L + b^{2-1/H}s) - \mathrm{E}\widetilde{Z}^{\mathcal{J}}(t_L)\widetilde{Z}^{\mathcal{J}}(t_L + b^{2-1/H}s))}{\operatorname{var}\widetilde{Z}(t_L) - \operatorname{var}\widetilde{Z}^{\mathcal{J}}(t_L)}$$
$$\quad \times (\mathrm{E}\widetilde{Z}(t_L)\widetilde{Z}(t_L + b^{2-1/H}t) - \mathrm{E}\widetilde{Z}^{\mathcal{J}}(t_L)\widetilde{Z}^{\mathcal{J}}(t_L + b^{2-1/H}t))$$
$$= \frac{\theta^2 \tilde{t}_L^{2H}}{2}((1+s/\tilde{t}_L)^{2H} + (1+t/\tilde{t}_L)^{2H} - (|t-s|/\tilde{t}_L)^{2H})$$
$$\quad - \theta^2 f(D/t_L, s/\tilde{t}_L, (t-s)/\tilde{t}_L)$$
$$\quad - \theta^2 \tilde{t}_L^{2H}\left[\frac{(v(s/\tilde{t}_L) - f(D/t_L, 0, s/\tilde{t}_L))(v(t/\tilde{t}_L) - f(D/t_L, 0, t/\tilde{t}_L))}{1 - f(D/t_L, 0, 0)}\right],$$

where $v(\cdot)$ and $f(\cdot, \cdot, \cdot)$ are as above. Straightforward algebraic manipulations then give

$$[\operatorname{cov}(\hat{Y}^b(s), \hat{Y}^b(t)|\mathcal{H}) - \operatorname{cov}(\hat{Y}(s), \hat{Y}(t))]\theta^{-2}\tilde{t}_L^{-2H}(1 - f(D/t_L, 0, 0))$$
$$= \frac{1 - f(D/t_L, 0, 0)}{2}$$
$$\quad \times \left[\left(1 + \frac{s}{\tilde{t}_L}\right)^{2H} - \left|\frac{s}{\tilde{t}_L}\right|^{2H}\right.$$
$$\quad \left. + (1+t/\tilde{t}_L)^{2H} - |t/\tilde{t}_L|^{2H} - 2f\left(\frac{D}{t_L}, \frac{s}{\tilde{t}_L}, (t-s)/\tilde{t}_L\right)\right]$$
$$\quad - \left[v\left(\frac{s}{\tilde{t}_L}\right) - f\left(\frac{D}{t_L}, 0, \frac{s}{\tilde{t}_L}\right)\right] \cdot \left[v(t/\tilde{t}_L) - f\left(\frac{D}{t_L}, 0, t/\tilde{t}_L\right)\right]$$
$$= (1 - f(D/t_L, 0, 0))\left[f(D/t_L, 0, t/\tilde{t}_L) + f\left(D/t_L, 0, \frac{s}{\tilde{t}_L}\right)\right.$$

$$- f(D/t_L, 0, 0) - f\left(D/t_L, \frac{s}{\tilde{t}_L}, (t-s)/\tilde{t}_L\right)\bigg]$$

$$- \left[f\left(D/t_L, 0, \frac{s}{\tilde{t}_L}\right) - f(D/t_L, 0, 0)\right][f(D/t_L, 0, t/\tilde{t}_L) - f(D/t_L, 0, 0)]$$

$$+ \left(v\left(\frac{s}{\tilde{t}_L}\right) - 1\right)[f(D/t_L, 0, t/\tilde{t}_L) - f(D/t_L, 0, 0)]$$

$$+ (v(t/\tilde{t}_L) - 1)\left[f\left(D/t_L, 0, \frac{s}{\tilde{t}_L}\right) - f(D/t_L, 0, 0)\right]$$

$$- \left(v\left(\frac{s}{\tilde{t}_L}\right) - 1\right)(v(t/\tilde{t}_L) - 1),$$

and the last four terms can be shown to be $O(b^{-2(1/H-1)}) = O(\tilde{t}_L^{-2})$ as $b \to \infty$, so that

$$[\text{cov}(\hat{Y}^b(s), \hat{Y}^b(t)|\mathcal{H}) - \text{cov}(\hat{Y}(s), \hat{Y}(t))]\theta^{-2}\tilde{t}_L^{-2H}$$
$$= f(D/t_L, 0, t/\tilde{t}_L) + f(D/t_L, 0, s/\tilde{t}_L) - f(D/t_L, 0, 0)$$
$$\quad - f(D/t_L, s/\tilde{t}_L, (t-s)/\tilde{t}_L) + O(\tilde{t}_L^{-2})$$
$$= H \int_0^{D/t_L} \{g(D/t_L, 1, -u)[(1 + u + t/\tilde{t}_L)^{2H-1}$$
$$\qquad\qquad + (1 + u + s/\tilde{t}_L)^{2H-1}]$$
$$\qquad - g(D/t_L, 1, -u)[(1+u)^{2H-1} + u^{2H-1}]$$
$$\qquad - g(D/t_L, 1 + s/\tilde{t}_L, -u)[(1 + u + t/\tilde{t}_L)^{2H-1} - u^{2H-1}]\} \, du$$
$$\quad + O(\tilde{t}_L^{-2})$$
$$= H \int_0^{D/t_L} \{g(D/t_L, 1, -u)[(1 + u + s/\tilde{t}_L)^{2H-1} - (1+u)^{2H-1}]$$
$$\qquad - (1 + u + t/\tilde{t}_L)^{2H-1}$$
$$\qquad\qquad \times [g(D/t_L, 1 + s/\tilde{t}_L, -u) - g(D/t_L, 1, -u)]\} \, du$$
$$\quad + O(\tilde{t}_L^{-1} D t_L^{-1}),$$

where the last step follows because

$$\left|\int_0^{D/t_L} u^{2H-1}[g(D/t_L, 1, -u) - g(D/t_L, 1 + s/\tilde{t}_L, -u)] \, du\right|$$
$$= C_H \int_0^{D/t_L} u^{H-1/2}(D/t_L - u)^{-H+1/2}$$

$$\times \int_{1\wedge(1+s/\tilde{t}_L)}^{1\vee(1+s/\tilde{t}_L)} \frac{y^{H-1/2}(y+D/t_L)^{H-1/2}}{y+u}\, dy\, du$$

$$\leq 2|s/\tilde{t}_L| \int_0^{D/t_L} u^{H-1/2}(D/t_L - u)^{-H+1/2}\, du$$

$$= 2|s/\tilde{t}_L|(D/t_L) \int_0^1 v^{H-1/2}(1-v)^{-H+1/2}\, dv$$

$$= O(\tilde{t}_L^{-1} D t_L^{-1})$$

as $b \to \infty$. Put

$$h(d, x, u) \triangleq \int_0^{1+x} \frac{y^{H-1/2}(y+d)^{H-1/2}}{y+u}\, dy,$$

and note that there exist constants $c_6$ and $c_7$ such that, for all small enough $d$ and all $0 < u, |x| < d$,

$$|(1+u+x)^{2H-1} - (1+u)^{2H-1} - (2H-1)(1+u)^{2H-2}x| \leq c_6 x^2,$$

$$|h(d, x, u) - h(d, 0, u) - x(1+d)^{H-1/2}(1+u)^{-1}| \leq c_7 x^2.$$

It follows that, for all large enough $b$,

$$|\text{cov}(\hat{Y}^b(s), \hat{Y}^b(t)|\mathcal{H}) - \text{cov}(\hat{Y}(s), \hat{Y}(t))|\theta^{-2}\tilde{t}_L^{-2H}$$

$$= H \int_0^{D/t_L} u^{-H+1/2}\left(\frac{D}{t_L} - u\right)^{-H+1/2}$$

$$\times \left|[(1+u+s/\tilde{t}_L)^{2H-1} - (1+u)^{2H-1}]h\left(\frac{D}{t_L}, 0, u\right)\right.$$

$$- (1+u+t/\tilde{t}_L)^{2H-1}$$

$$\left.\times [h(D/t_L, s/\tilde{t}_L, -u) - h(D/t_L, 0, -u)]\right|\, du$$

$$+ O(\tilde{t}_L^{-1} D t_L^{-1}),$$

$$= H \int_0^{D/t_L} u^{-H+1/2}(D/t_L - u)^{-H+1/2}$$

$$\times |(2H-1)(1+u)^{2H-2}h(D/t_L, 0, u)(s/\tilde{t}_L)$$

$$- (1+u+t/\tilde{t}_L)^{2H-1}(1+D/t_L)^{H-1/2}(1+u)^{-1}(s/\tilde{t}_L)|\, du$$

$$+ O((s/\tilde{t}_L)^2 (D/t_L)^{2(1-H)}) + O(\tilde{t}_L^{-1} D t_L^{-1}).$$

But there exist constants $c_8$, $c_9$ such that, for all small enough $d$ and all $0 < u$, $|y| < d$,

$$(1+d)^{H-1/2} - 1 \leq c_8 d \quad \text{and} \quad |(1+u+y)^{2H-1} - (1+u)^{2H-1}| \leq c_9|y|,$$

so that, for all large enough $b$,

$$|\operatorname{cov}(\hat{Y}^b(s), \hat{Y}^b(t)|\mathcal{H}) - \operatorname{cov}(\hat{Y}(s), \hat{Y}(t))|\theta^{-2}\tilde{t}_L^{-2H}$$

$$\leq H(s/\tilde{t}_L) \int_0^{D/t_L} u^{-H+1/2} \left(\frac{D}{t_L} - u\right)^{-H+1/2} (1+u)^{2H-2}$$

$$\times \left|(2H-1)h\left(\frac{D}{t_L}, 0, u\right) - 1\right| du$$

$$+ O((s/\tilde{t}_L)(D/t_L)^{3-2H}) + O((s/\tilde{t}_L)^2(D/t_L)^{2(1-H)}) + O(\tilde{t}_L^{-1}Dt_L^{-1}).$$

But, for all small enough $d > 0$, $|(2H-1)h(d,0,u) - 1| \leq 5(1-H)^{-1}d^{2H-1}$, so that for an appropriate constant $c_{10}$ and all large enough $b$,

$$|\operatorname{cov}(\hat{Y}^b(s), \hat{Y}^b(t)|\mathcal{H}) - \operatorname{cov}(\hat{Y}(s), \hat{Y}(t))|\theta^{-2}\tilde{t}_L^{-2H}$$

$$\leq c_{10}(s/\tilde{t}_L)\left(\frac{D}{t_L}\right)^{2H-1} \int_0^{D/t_L} u^{-H+1/2}\left(\frac{D}{t_L} - u\right)^{-H+1/2} (1+u)^{2H-2} du$$

$$+ O(\tilde{t}_L^{-1}Dt_L^{-1})$$

$$= c_{10}(s/\tilde{t}_L)(D/t_L)^{2H-1}O((D/t_L)^{2(1-H)}) + O(\tilde{t}_L^{-1}Dt_L^{-1})$$

$$= O(\tilde{t}_L^{-1}(D/t_L)),$$

and we conclude that

$$|\operatorname{cov}(\hat{Y}^b(s), \hat{Y}^b(t)|\mathcal{H}) - \operatorname{cov}(\hat{Y}(s), \hat{Y}(t))| = O(\tilde{t}_L^{2H-1}(D/t_L))$$

$$= O(b^{-(1-H)^2/H+r})$$

$$= o(b^{-\varepsilon_2})$$

for any $0 < \varepsilon_2 < -(1-H)^2/H + r$. $\square$

LEMMA 4. *Assume (5) and (4) hold. Given $M > 0$, for all large enough $b$ the function $u \mapsto E(\hat{Y}^b(u)|\mathcal{H})$ is differentiable on $[-M, M]$ and there exists a constant $\widetilde{C}_M$ such that, on $A(b)$,*

$$\left|\frac{\partial E(\hat{Y}^b(u)|\mathcal{H})}{\partial u}\right| \leq \widetilde{C}_M, \tag{22}$$

$u \in [-M, M]$.

PROOF. Recall that, with $f$ and $v$ as in the proof of Lemma 3,

$$\mathrm{E}(\hat{Y}^b(u)|\mathcal{H}) = b^{1-2H}\frac{f(D/t_L, 0, u/\tilde{t}_L) - v(u/\tilde{t}_L)}{1 - f(D/t_L, 0, 0)}\widetilde{Z}_{\mathcal{G}}^{\mathcal{J}}(t_L)$$

$$+ b^{1-2H}\widetilde{Z}_{\mathcal{G}}^{\mathcal{J}}(t_L + b^{2-1/H}u) \tag{23}$$

$$+ [\mathrm{E}(\hat{Y}^b(u)|\mathcal{G}) - \mathrm{E}(\hat{Y}^b(u)|\mathcal{K})] + \mathrm{E}(\hat{Y}^b(u)|\mathcal{K}).$$

The first term on the RHS in (23) is differentiable (as a function of $u$), and

$$\frac{\partial}{\partial u}\left(b^{-2H+1}\frac{f(D/t_L, 0, u/\tilde{t}_L) - v(u/\tilde{t}_L)}{1 - f(D/t_L, 0, 0)}\widetilde{Z}_{\mathcal{G}}^{\mathcal{J}}(t_L)\right)$$

$$= \left[\frac{(\partial f/\partial x_3)(D/t_L, 0, u/\tilde{t}_L) - v'(u/\tilde{t}_L)}{1 - f(D/t_L, 0, 0)}\right][b^{1-2H}\widetilde{Z}_{\mathcal{G}}^{\mathcal{J}}(t_L)\tilde{t}_L^{-1}].$$

Here, the first factor on the RHS is bounded over $[-M, M]$ for large enough $b$, while the second factor does not depend on $u$ and is $o(1)$ as $b \to \infty$, as can be seen from (19) in the proof of Lemma 2. Hence, there exists a constant $\tilde{C}_{M1}$ such that, on $A(b)$,

$$\left|\frac{\partial}{\partial u}\left(b^{-2H+1}\frac{v(u/\tilde{t}_L) - f(D/t_L, 0, u/\tilde{t}_L)}{1 - f(D/t_L, 0, 0)}\widetilde{Z}_{\mathcal{G}}^{\mathcal{J}}(t_L)\right)\right| \leq \tilde{C}_{M1},$$

$u \in [-M, M]$.

For the second term on the RHS in (23), the integral representation used in the proof of Lemma 2 can be used to obtain

$$b^{1-2H}\frac{\partial}{\partial u}\widetilde{Z}_{\mathcal{G}}^{\mathcal{J}}(t_L + b^{2-1/H}u)$$

$$= b^{1-2H}\frac{C_H}{t_L}\int_0^{D/t_L} x^{1/2-H}\left(\frac{D}{t_L} - x\right)^{1/2-H}$$

$$\times \frac{(1 + u/\tilde{t}_L)^{H-1/2}(1 + D/t_L + u/\tilde{t}_L)^{H-1/2}}{(1 + x + u/\tilde{t}_L)}$$

$$\times \left[\frac{H - 1/2}{x} - \frac{H - 1/2}{D/t_L - x} - \frac{1}{1 + x + u/\tilde{t}_L}\right]\widetilde{Z}_{\mathcal{G}}(x\tilde{t}_L)\,dx.$$

Using that, on $A(b)$, $|\widetilde{Z}_{\mathcal{G}}(-w)| \leq \max(w^H b^\delta, w^{H-\eta_1})$ for $0 \leq w \leq D/2$ and $|\widetilde{Z}_{\mathcal{G}}(-w)| \leq \max((D-w)^H b^\delta, (D-w)^{H-\eta_1})$ for $D/2 \leq w \leq D$, one can conclude that the RHS on the above equation is bounded, so that there exists $\tilde{C}_{M2}$ such that, on $A(b)$,

$$b^{1-2H}\left|\frac{\partial}{\partial u}\widetilde{Z}_{\mathcal{G}}^{\mathcal{J}}(t_L - b^{2-1/H}u)\right| \leq \tilde{C}_{M2},$$

$u \in [-M, M]$.

For the third term on the RHS in (23), recall

$$\mathrm{E}(\hat{Y}^b(u)|\mathcal{G}) - \mathrm{E}(\hat{Y}^b(u)|\mathcal{K})$$

$$= \theta b^{1-2H} \cdot \frac{Z_{\mathcal{K}}(-t_R)}{4(1 + D/t_L)^{2H} - (1 + (1 + D/t_L)^{2H} - (D/t_L)^{2H})^2}$$

$$\times \left[2\left|\frac{u}{\tilde{t}_L}\right|^{2H} + 2\left(1 + \frac{D}{t_L}\right)^{2H} - 2\left(1 + \frac{D}{t_L} + \frac{u}{\tilde{t}_L}\right)^{2H}\right.$$

$$-\left(\left|\frac{u}{\tilde{t}_L}\right|^{2H}+1-\left(1+\frac{u}{\tilde{t}_L}\right)^{2H}\right)\left(1+\left(1+\frac{D}{t_L}\right)^{2H}-\left(\frac{D}{t_L}\right)^{2H}\right)\right],$$

whence

$$\left|\frac{\partial}{\partial u}(\mathrm{E}(\hat{Y}^b(u)|\mathcal{G})-\mathrm{E}(\hat{Y}^b(u)|\mathcal{K}))\right|$$

$$= 2H\tilde{t}_L^{-1}\theta b^{1-2H}\left|\frac{Z_\mathcal{K}(-t_R)}{4(1+D/t_L)^{2H}-(1+(1+D/t_L)^{2H}-(D/t_L)^{2H})^2}\right|$$

$$\times \left|2\,\mathrm{sign}(u)\left|\frac{u}{\tilde{t}_L}\right|^{2H-1}-2\left(1+\frac{D}{t_L}+\frac{u}{\tilde{t}_L}\right)^{2H-1}\right.$$

$$-\left(1+\left(1+\frac{D}{t_L}\right)^{2H}-\left(\frac{D}{t_L}\right)^{2H}\right)$$

$$\left.\times\left(\mathrm{sign}(u)\left|\frac{u}{\tilde{t}_L}\right|^{2H-1}-\left(1+\frac{u}{\tilde{t}_L}\right)^{2H-1}\right)\right|$$

$$\leq (1+M^{2H-1})O(b^{-\varepsilon_2}),$$

where the term $O(b^{-\varepsilon_2})$ does not depend on $u$, and $\varepsilon_2$ is as in Step 2. Hence, there exists a constant $\tilde{C}_{M3}$ such that, on $A(b)$,

$$\left|\frac{\partial}{\partial u}(\mathrm{E}(\hat{Y}^b(u)|\mathcal{G})-\mathrm{E}(\hat{Y}^b(u)|\mathcal{K}))\right|\leq \tilde{C}_{M3},$$

$u\in[-M,M]$.

For the last term on the RHS in (23), one can differentiate the expression for $\mathrm{E}(\hat{Y}^b(u)|\mathcal{K})$ given in (3) (with $u$ in place of $s$) to obtain

$$\frac{\partial}{\partial u}\mathrm{E}(\hat{Y}^b(u)|\mathcal{K}) = -2H\kappa_1|u|^{2H-1}\mathrm{sign}(u)(t^*/t_L)^{2H}$$

$$+\frac{Hb^{2(1-H)}}{(1-H)\tilde{t}_L}(1-(1+u/\tilde{t}_L)^{2H-1})+\mu\tilde{t}_L^{-1}b^{1-H+r}$$

$$+ Hb^{1-2H}\tilde{t}_L^{-1}[X(-t_L)+b+\mu b^{H+r}]$$

$$\times (|u/\tilde{t}_L|^{2H-1}\mathrm{sign}(u)-(1+u/\tilde{t}_L)^{2H-1}),$$

and since $|X(-t_L)+b|\leq b^\gamma < b^{H+r}$ on $A(b)$, we conclude there exists a constant $\tilde{C}_{M4}$ such that, on $A(b)$,

$$\left|\frac{\partial}{\partial u}\mathrm{E}(\hat{Y}^b(u)|\mathcal{K})\right|\leq \tilde{C}_{M4},$$

$u\in[-M,M]$. Taking $\tilde{C}_M = \tilde{C}_{M1}+\tilde{C}_{M2}+\tilde{C}_{M3}+\tilde{C}_{M4}$, the result follows.
□

LEMMA 5. *If $\gamma$, $\eta_1$, $r$ and $\delta$ satisfy (15)–(17), then*

$$P\{A(b) \triangle \{T^b < \infty\}\} = o(P\{T^b < \infty\}).$$

PROOF. We show that all the terms on the RHS in (14) are $o(P(T^b < \infty))$ as $b \to \infty$.

For the first term, note $P(A_5) \leq P(T^b < t_L) = o(P(T^b < \infty))$ by Theorem 3.

For the second term, choose any $\varepsilon$ such that $0 < \varepsilon < (2H-1)(1-H-r)/8$, and put $\nu_1 \triangleq 1 - H + H^2 + rH + \varepsilon$. Note that

$$P(A_0 \cap \{X(-t_L) > -b + b^{\nu_1}\}) \leq P\Big(\sup_{0 \leq s \leq 2b^{H+r}} Z(-t_L) - Z(-t_L - s) > b^{\nu_1}/\theta\Big),$$

and it then follows from Borell's inequality (Lemma 7 below) that there exist constants $c_1$ and $c_2$ such that

$$P(A_0 \cap \{X(-t_L) > -b + b^{\nu_1}\}) \leq c_1 \exp(-c_2 b^{2\nu_1}/b^{2H(H+r)}) = o(P(T^b < \infty)),$$

where the last follows since $2\nu_1 - 2H(r+H) > 2(1-H)$ and

$$P(T^b < \infty) \sim q(b) \exp(-(\kappa_2/(2(1-H)))b^{2(1-H)}),$$

where $q(\cdot)$ is a ratio of polynomials, by a result of Hüsler and Piterbarg [15], Corollary 2.

Let $\nu_1$ and $\varepsilon$ be as above, and for $k \geq 1$ put

$$\nu_{k+1} = (H + r - (2H-1)(1-H-r))/2 + (\nu_k \vee (H+r))/2 + \varepsilon.$$

Observe that $\nu_1 > r + H$, and it follows from Lemma 6 that

$$P(A_0 \cap \{X(-t_L) > -b + b^{\nu_2}\}) = o(P(T^b < \infty))$$

as $b \to \infty$. Noting that $\nu_k$ is nonincreasing in $k$, an inductive application of Lemma 6 shows that whenever $\nu_k > H + r$, then

$$P(A_0 \cap \{X(-t_L) > -b + b^{\nu_{k+1}}\}) = o(P(T^b < \infty)).$$

Also, note that if $k \geq 2$ and $\nu_k > H + r$, then

$$\nu_k - \nu_{k+1} = (\nu_{k-1} - \nu_k)/2.$$

It then follows by an inductive argument that, if $\nu_k > H + r$, then

$$\nu_{k+1} = H + r - (2H-1)(1-H-r) + 2^{-k}H(1-H-r) + 2(1-2^{-k+1})\varepsilon.$$

Let $k^* \triangleq \inf\{k : \nu_k \leq H + r\}$ and note $1 < k^* < \infty$. Also, note that $\nu_{k^*} \geq H + r - (2H-1)(1-H-r)/2$. Then, putting $\nu^* \triangleq \nu_{k^*+1} = H + r - (H - 1/2)(1-H-r) + \varepsilon$, it follows by Lemma 6 that

$$P(A_0 \cap \{X(-t_L) > -b + b^{\nu^*}\}) = o(P(T^b < \infty)).$$

If $\gamma$ satisfies (15) and $\varepsilon$ is small enough, then $\nu^* < \gamma$, so that the second term on the RHS in (14) satisfies

$$P(A_0 \cap A_6) = P(A_0 \cap \{X(-t_L) > -b + b^\gamma\})$$
$$\leq P(A_0 \cap \{X(-t_L) > -b + b^{\nu^*}\}) = o(P(T^b < \infty)).$$

For the third term on the RHS in (14), note that

$$P(A_1 \cap A_2^C)$$
$$= EI_{A_1} P(|Z_\mathcal{K}(-t_R)| > b^{\eta_3}|\mathcal{K})$$
$$\leq P(X(-t_L) < -b + b^\gamma) \cdot P(2^H b^{H(H+r)}|N(0,1)| > b^{\eta_3})$$

But it follows from the Hüsler–Piterbarg asymptotic (13) and equation (26) in the Proof of Lemma 6 that

$$P(X(-t_L) < -b + b^\gamma) = P(T^b < \infty) q(b) \exp(b^{2(1-H)} O(b^{-1+H+r})),$$

while

$$P(2^H b^{H(H+r)}|N(0,1)| > b^{\eta_3}) = O(b^{-\eta_3 + H(H+r)} \exp(-b^{2(\eta_3 - H^2 - rH)}/2^{2H-1}))$$

in view of (16), it then follows that $P(A_1 \cap A_2^C) = o(P(T^b < \infty))$.

To deal with the fourth term on the RHS in (14), put $N \triangleq \lceil b^{\eta_6} \rceil$, $h \triangleq b^{-2(1-H)/\eta_1}$, and $\Delta \triangleq (t^* - t_L - h)/N$, where $\eta_6 > 0$. Note that

(24)
$$P(A_1 \cap A_3^C) \leq \sum_{n=0}^N P(A_1 \cap E_n) + \sum_{n=0}^N P(A_1 \cap B_n)$$
$$+ \sum_{n=1}^\infty P(A_1 \cap C_n) + \sum_{n=2}^\infty P(A_1 \cap D_n),$$

where

$E_n \triangleq \{|Z_\mathcal{G}(-t_L - h - n\Delta)| > (1/2) b^\delta (h + n\Delta)^H\}$,

$B_n \triangleq \left\{\sup_{0 \leq s \leq \Delta} |Z_\mathcal{G}(-t_L - h - n\Delta) - Z_\mathcal{G}(-t_L - h - n\Delta - s)| > (1/2) b^\delta h^H\right\}$,

$C_n \triangleq \{|Z_\mathcal{G}(-t_L - h/n^2)| > (1/2)(h/n^2)^{H-\eta_1}\}$,

$D_n \triangleq \left\{\sup_{hn^{-2} \leq s \leq h(n-1)^{-2}} |Z_\mathcal{G}(-t_L - s) - Z_\mathcal{G}(-t_L - h/n^2)| > (1/2)(h/n^2)^{H-\eta_1}\right\}$.

Note

(25)
$$P(A_1 \cap E_n)$$
$$= EI_{A_1} E(I_{E_n}|\mathcal{G})$$
$$\leq EI_{A_1} P((h + n\Delta)^H |N(0,1)| > (1/2)(h + n\Delta)^H b^\delta)$$
$$\leq P(X(-t_L) < -b + b^\gamma) O(b^{-\delta} \exp(-b^{2\delta}/8))$$
$$= P(T^b < \infty) q(b) \exp(b^{2(1-H)} O(b^{-1+H+r})) O(b^{-\delta} \exp(-b^{2\delta}/8)),$$

where $q(b)$ is a ratio of polynomials and the last step follows from the Hüsler–Piterbarg asymptotic (13) and equation (26) in the Proof of Lemma 6. Since $\delta > (1-H+r)/2$, it follows that there exists a constant $c_3 > 0$ independent of $n$ such that

$$P(A_1 \cap E_n) \leq P(T^b < \infty) q(b) b^{-\delta} O(\exp(O(b^{1-H+r}) - b^{2\delta}/8)),$$
$$\leq P(T^b < \infty) \exp(-c_3 b^{2\delta}).$$

Hence,

$$\sum_{n=0}^{N} P(A_1 \cap E_n) \leq P(T^b < \infty) 2 b^{\eta_6} \exp(-c_3 b^{2\delta}) = o(P(T^b < \infty))$$

as $b \to \infty$.

Next, note that on $A_1$,

$$|Z^{\mathcal{G}}(-t_L - h - n\Delta + s) - Z^{\mathcal{G}}(-t_L - h - n\Delta)| \leq c_4 \Delta$$

for some constant $c_4 > 0$, all $0 \leq s \leq \Delta$ and all $0 \leq n \leq N$. Thus, provided $\eta_6$ is large enough so that $\Delta^H = o(b^\delta h^H)$ as $b \to \infty$, it follows that for large enough $b$

$$P(A_1 \cap B_n)$$
$$\leq P\left(\sup_{0 \leq s \leq \Delta} |Z_{\mathcal{G}}(-t_L - h - n\Delta - s) - Z_{\mathcal{G}}(-t_L - h - n\Delta)| > (1/4) b^\delta h^H\right)$$
$$\leq c_5 \exp(-c_6 b^{2\delta} h^{2H} \Delta^{-2H}),$$

for appropriate positive constants $c_5$ and $c_6$, by Borell's inequality. Hence, by choosing $\eta_6$ large enough we can ensure

$$P(A_1 \cap B_n) \leq \exp(-b^{4(1-H)}),$$

and it then follows that

$$\sum_{n=0}^{N} P(A_1 \cap B_n) \leq 2 b^{\eta_6} \exp(-b^{4(1-H)}) = o(P(T^b < \infty))$$

as $b \to \infty$.

A similar argument as the one used to analyze $P(A_1 \cap E_n)$ shows that there exist positive constants $c_1, c_2, c_3$ such that

$$P(A_1 \cap C_n) \leq c_1 \exp(-c_2 h^{-2\eta_1}) \exp(-c_3 n^{4\eta_1}).$$

Also, a similar argument as the one used to analyze $P(A_1 \cap B_n)$ shows that there exist positive constants $c_4, c_5, c_6$ such that

$$P(A_1 \cap D_n) \leq c_4 \exp(-c_5 h^{-2\eta_1}) \exp(-c_6 n^{4\eta_1}).$$

It then follows that the last two summations in (24) are of order $O(\exp(-(c_2 \wedge c_5)h^{-2\eta_1})) = o(\mathrm{P}(T_b < \infty))$, so that $\mathrm{P}(A_1 \cap A_3^C) = o(\mathrm{P}(T^b < \infty))$ as $b \to \infty$. The proof that $\mathrm{P}(A_1 \cap A_4^C) = o(\mathrm{P}(T^b < \infty))$ uses the same type of argument, and is omitted for brevity.

Hence, all terms on the RHS in (14) are $o(\mathrm{P}(T^b < \infty))$, what had to be shown. $\square$

LEMMA 6. *Let $r$ be the small number chosen in the definition of $A_0$. If $1 > \nu > H + r - (2H-1)(1-H-r)/2$ and $\mathrm{P}(A_0 \cap \{X(-t_L) > -b + b^\nu\}) = o(\mathrm{P}(T^b < \infty))$, then for any $\tilde{\nu} > (H + r - (2H-1)(1-H-r))/2 + (\nu \vee (H+r))/2$ it also holds that*

$$\mathrm{P}(A_0(b) \cap \{X(-t_L) > -b + b^{\tilde{\nu}}\}) = o(\mathrm{P}(T^b < \infty))$$

*as $b \to \infty$.*

PROOF. Observe that

$$\mathrm{P}(A_0(b) \cap \{X(-t_L) > -b + b^{\tilde{\nu}}\})$$
$$= \mathrm{P}(A_0(b) \cap \{X(-t_L) \in [-b + b^{\tilde{\nu}}, -b + b^\nu]\}) + o(\mathrm{P}(T^b < \infty))$$
$$= \mathrm{E}[I(X(-t_L) \in [-b + b^{\tilde{\nu}}, b - b^\nu]) \cdot \mathrm{P}(A_0(b)|\mathcal{K})] + o(\mathrm{P}(T^b < \infty)).$$

But on $\{X(-t_L) \in [-b + b^{\tilde{\nu}}, -b + b^\nu]\}$, $\mathrm{E}(X(-t_L) - X(-t_L - s)|\mathcal{K}) = o(1)$ as $b \to \infty$. In particular, for all large enough $b$,

$$I(X(-t_L) \in [-b + b^{\tilde{\nu}}, -b + b^\nu]) \cdot I_{A_0(b)}$$
$$\leq I\left(\sup_{0 \leq s \leq 2b^{H+r}} -Z_\mathcal{K}(-t_L - s) > b^{\tilde{\nu}}/(2\theta)\right).$$

By Borell's inequality, we can find positive constants $c_1$ and $c_2$ such that

$$\mathrm{P}\left(\sup_{0 \leq s \leq 2b^{H+r}} -Z_\mathcal{K}(-t_L - s) > b^{\tilde{\nu}}/(2\theta)\Big|\mathcal{K}\right) \leq c_1 \exp(-c_2 b^{2\tilde{\nu}}/b^{2H(H+r)})$$

on $\{X(t_L) \in [-b + b^{\tilde{\nu}}, -b + b^\nu]\}$. It then follows that

$$\mathrm{P}(A_0(b) \cap \{X(-t_L) > -b + b^{\tilde{\nu}}\})$$
$$\leq \mathrm{P}(X(-t_L) < -b + b^\nu) \cdot c_1 \exp(-c_2 b^{2\tilde{\nu} - 2H(H+r)}) + o(\mathrm{P}(T^b < \infty)).$$

But

(26)
$$\mathrm{P}(X(-t_L) < -b + b^\nu)$$
$$\sim \theta(1-H)b^{-1}\exp\left\{-\frac{\kappa_2}{2(1-H)}b^{2(1-H)}(1-h(b))\right\},$$

where

$$h(b) \triangleq 1 - \frac{(1 - (\mu/(1-H))b^{H+r-1} - b^{\nu-1}/(1-H))^2}{(1 - (\mu(1-H)/H)b^{H+r-1})^{2H}} = O(b^{-1+\nu \vee (H+r)})$$

as $b \to \infty$. Hence, with $q(\cdot)$ a ratio of polynomials,

$$\begin{aligned}
&\mathrm{P}(A_0(b) \cap \{X(-t_L) > -b + b^{\tilde\nu}\}) \\
&\quad \leq \mathrm{P}(T^b < \infty)\left\{q(b)\exp\left(\frac{\kappa_2 h(b)}{2(1-H)}b^{2(1-H)}\right)\exp(-c_2 b^{2\tilde\nu - 2H(H+r)}) + o(1)\right\} \\
&\quad \leq \mathrm{P}(T^b < \infty)\{q(b)\exp(O(b^{1-2H+\nu \vee (H+r)})) - c_2 b^{2\tilde\nu - 2H(H+r)}) + o(1)\} \\
&\quad = o(\mathrm{P}(T^b < \infty))
\end{aligned}$$

as $b \to \infty$, since $1 - 2H + \nu \vee (H+r) < 2\tilde\nu - 2H(H+r)$. $\square$

The following result is known as Borell's inequality in the literature; see Adler [1], Theorem 2.1. We quote it here since we use it at many points in the proofs of Theorems 1 and 3.

LEMMA 7 (Borell's inequality). *Let $\{V(t) : t \in T\}$ be a centered Gaussian process with sample paths bounded a.s. Let $V^M = \sup_{t \in T} V(t)$. Then $\mathrm{E}V^M < \infty$ and for all $\lambda > 0$*

$$\mathrm{P}(|V^M - \mathrm{E}V^M| > \lambda) \leq 2e^{-(1/2)\lambda^2/\sigma_T^2},$$

*where $\sigma_T^2 = \sup_{t \in T} \mathrm{var}\, V(t)$. As an immediate consequence, for $\lambda > \mathrm{E}V^M$,*

$$\mathrm{P}(V^M > \lambda) \leq 2e^{-(1/2)(\lambda - \mathrm{E}V^M)^2/\sigma_T^2}.$$

*In particular, if $Z$ is standard fBM and $C \triangleq \mathrm{E}\sup_{0 \leq s \leq 1} Z(s)$, then for $\lambda > t^H C$*

$$\mathrm{P}\Big(\sup_{s \in [0,t]} Z(s) > \lambda\Big) \leq 2e^{-(1/2)(\lambda - t^H C)^2/t^{2H}}.$$

Table 1 summarizes selected notation which is used repeatedly at different points of the argument.

**Acknowledgments.** The authors gratefully acknowledge the financial support of Alcatel Research and Innovation (now Alcatel–Lucent) over the course of this research. In addition, we thank Peter Rabinovitch of Alcatel–Lucent for his many valuable suggestions concerning this paper, and two anonymous referees whose comments greatly improved the presentation.

Table 1

$$
\begin{aligned}
&t_L = t^* - b^{H+r}; \quad t_R = t^* + b^{H+r}; \quad \tilde{t}_L = b^{-2+1/H}t_L; \quad D = 2b^{H+r} \\
&\hat{Y}^b(t) = b^{-2H+1}X(b^{2-1/H}t) \\
&\widetilde{Z}(t) = Z(-t_L + t) - Z(-t_L) \\
&\mathcal{K} = \sigma(Z(-t_L)) = \sigma(\widetilde{Z}(t_L)) \\
&\mathcal{G} = \sigma(Z(-t_L), Z(-t_R)) = \sigma(\widetilde{Z}(t_L), \widetilde{Z}(-D)) \\
&\mathcal{H} = \sigma(Z(s)\colon -t_R < s < -t_L) = \sigma(\widetilde{Z}(t_L), (\widetilde{Z}(s)\colon -D < s < 0)) \\
&\mathcal{J} = \sigma(\widetilde{Z}(s)\colon -D < s < 0) \\
&A_0(b) = \{\exists s \in [-t_R, -t_L]\colon X(s) = -b\} \\
&A_1(b) = \{X(-t_L) \in [-b, -b+b^{\gamma}]\} \\
&A_2(b) = \{|Z_{\mathcal{K}}(-t_R)| \leq b^{\eta_3}\} \\
&A_3(b) = \{|Z_{\mathcal{G}}(-t_L - s)| \leq \max(s^H b^{\delta}, s^{H-\eta_1}), 0 \leq s \leq b^{H+r}\} \\
&A_4(b) = \{|Z_{\mathcal{G}}(-t_R + s)| \leq \max(s^H b^{\delta}, s^{H-\eta_1}), 0 \leq s \leq b^{H+r}\} \\
&A(b) = A_0(b) \cap A_1(b) \cap A_2(b) \cap A_3(b) \cap A_4(b)
\end{aligned}
$$

DEPARTMENT OF MANAGEMENT SCIENCE
UNIVERSITY OF MIAMI
CORAL GABLES, FLORIDA 33146-6544
USA
E-MAIL: h.awad@miami.edu

DEPARTMENT OF MANAGEMENT SCIENCE
AND ENGINEERING
STANFORD UNIVERSITY
STANFORD, CALIFORNIA 94305-4026
USA
E-MAIL: glynn@stanford.edu